\documentclass{amsart}

\usepackage{amsfonts}
\usepackage{amsmath}
\usepackage{amsthm}
\usepackage{amssymb}
\usepackage[english]{babel}
\usepackage{graphics}
\usepackage{latexsym}
\usepackage{longtable} 
\usepackage{mathrsfs} 
\usepackage{graphicx}
\usepackage{wrapfig} 
\usepackage[pdftex,colorlinks=true,linkcolor=blue,urlcolor=blue,citecolor=blue]{hyperref}

\newtheorem{theorem}{Theorem}
\newtheorem{corollary}[theorem]{Corollary}
\newtheorem{lemma}[theorem]{Lemma}
\newtheorem{proposition}[theorem]{Proposition}

\newtheorem{claim}[theorem]{Claim}
\newtheorem{example}[theorem]{Example}

\theoremstyle{definition}
\newtheorem{definition}[theorem]{Definition}
\newtheorem{remark}[theorem]{Remark}

\newcommand{\mL}{\mathcal{L}}
\newcommand{\mH}{\mathcal{H}}
\newcommand{\mF}{\mathcal{F}}

\newcommand{\mK}{\mathbb{K}}
\newcommand{\mW}{\mathbb{W}}

\newcommand{\mP}{\mathscr{P}}
\newcommand{\mM}{\mathcal{M}}

\newcommand{\mY}{\mathscr{Y}}
\newcommand{\mD}{\mathcal{D}}

\newcommand{\A}{\textbf{A}}
\newcommand{\R}{\mathbb{R}}

\newcommand{\N}{\mathbb{N}}

\newcommand{\mB}{\mathbb{B}}
\newcommand{\X}{\textbf{X}}
\newcommand{\Y}{\textbf{Y}}
\newcommand{\Z}{\textbf{Z}}

\newcommand{\noi}{\noindent}
\newcommand{\ms}{\medskip}

\newcommand{\al}{\alpha}
\newcommand{\be}{\beta}
\newcommand{\ga}{\gamma}
\newcommand{\Ga}{\Gamma}
\newcommand{\de}{\delta}
\newcommand{\De}{\Delta}
\newcommand{\e}{\varepsilon}
\newcommand{\si}{\sigma}
\newcommand{\Si}{\Sigma}
\newcommand{\la}{\lambda}
\newcommand{\La}{\Lambda}
\newcommand{\ka}{\kappa}
\newcommand{\Om}{\Omega}
\newcommand{\om}{\omega}

\newcommand{\Th}{\Theta}

\newcommand{\D}{\mathrm{D}} 
 
\newcommand{\weak }{\hspace{1pt} -\!\!\!\!-\!\!\!\rightharpoonup}
\newcommand{\weakstar }{ \overset{\hspace{1pt} *_{\phantom{|}}}{{\smash{\weak }}\hspace{1pt} } }

\newcommand{\larrow}{\longrightarrow}
\newcommand{\ot}{\otimes}

\newcommand{\lmapsto}{\longmapsto}
\newcommand{\ri}{\rightarrow}
\newcommand{\LL}{\mathsf{L}}
\newcommand{\p}{\partial}
\newcommand{\sub}{\subseteq}
\newcommand{\set}{\setminus}
\newcommand{\by}{\times}

\newcommand{\Lip}{\mathrm{Lip}}

\newcommand{\sgn}{\mathrm{sgn}}
\newcommand{\ess}{\mathrm{ess}}
\newcommand{\diam}{\mathrm{diam}}
\newcommand{\dist}{\mathrm{dist}}

\newcommand{\spn}{\mathrm{span}}
\newcommand{\supp}{\mathrm{supp}}

\newcommand{\bt}{\begin{theorem}}\newcommand{\et}{\end{theorem}}
\newcommand{\bd}{\begin{definition}}\newcommand{\ed}{\end{definition}}
\newcommand{\bl}{\begin{lemma}}\newcommand{\el}{\end{lemma}}
\newcommand{\beq}{\begin{equation}}\newcommand{\eeq}{\end{equation}}
\newcommand{\bc}{\begin{claim}}\newcommand{\ec}{\end{claim}}
\newcommand{\bex}{\begin{example}}\newcommand{\eex}{\end{example}}
\newcommand{\bcor}{\begin{corollary}}\newcommand{\ecor}{\end{corollary}}
\newcommand{\bp}{\begin{proof}}\newcommand{\ep}{\end{proof}}

\newcommand{\BPL}{\medskip \noindent \textbf{Proof of Lemma} }
\newcommand{\BPC}{\medskip \noindent \textbf{Proof of Claim} }

\newcommand{\BPP}{\medskip \noindent \textbf{Proof of Proposition} }
\newcommand{\BPT}{\medskip \noindent \textbf{Proof of Theorem} }

\numberwithin{equation}{section}

\begin{document}

\title[Generalised Solutions for Fully Nonlinear Systems and Existence]{Generalised Solutions for Fully Nonlinear PDE Systems and Existence-Uniqueness Theorems}

\author{Nikos Katzourakis}
\address{Department of Mathematics and Statistics, University of Reading, Whiteknights, PO Box 220, Reading RG6 6AX, Berkshire, UK}
\email{n.katzourakis@reading.ac.uk}

    \thanks{\!\!\!\!\!\!\texttt{The author has been partially financially supported by the EPSRC grant EP/N017412/1}}



\date{}



\begin{abstract} We introduce a new theory of generalised solutions which applies to fully nonlinear PDE systems of any order and allows for merely measurable maps as solutions. This approach bypasses the standard problems arising by the application of Distributions to PDEs and is not based on either integration by parts or on the maximum principle. Instead, our starting point builds on the probabilistic representation of derivatives via limits of difference quotients in the Young measures over a toric compactification of the space of jets. After developing some basic theory, as a first application we consider the Dirichlet problem and we prove existence-uniqueness-partial regularity of solutions to fully nonlinear degenerate elliptic 2nd order systems and also existence of solutions to the $\infty$-Laplace system of vectorial Calculus of Variations in $L^\infty$.
\end{abstract}

\maketitle

\tableofcontents

\section{Introduction} \label{section1}

It is well known that PDEs, either linear or nonlinear, in general do not possess classical solutions, in the sense that not all derivatives that appear in the equation may actually exist. The standard approach to this problem consists of looking for appropriately defined \emph{generalised solutions} for which at least existence can be proved given certain boundary conditions. Subsequent considerations typically include uniqueness, regularity, qualitative properties and numerics. This approach has been enormously successful but unfortunately only PDEs with fairly special structure have been considered so far. A standing idea in this regard consists of using integration-by-parts in order to interpret derivatives ``weakly" by ``passing them to test functions". This duality method which dates back to the 1930s (\cite{S1, S2, So}) is basically restricted to divergence structure equations and systems. A more recent approach discovered in the 1980s is that of viscosity solutions (\cite{CL}) and builds on the maximum principle as a device to define ``weak" solutions. Although it applies mostly to the scalar case, it has been hugely successful since it includes fully nonlinear single equations. 

In this paper we introduce a new theory of generalised solutions which applies to nonlinear PDE systems of any order. Our approach allows for  \emph{merely measurable} maps to be rigorously interpreted and studied as solutions of systems which are possibly fully nonlinear and with discontinuous coefficients. More precisely, let $p,n,N,M \in \N$, let also $\Om \sub \R^n$ be an open set and 
\beq \label{1.1}
\ \mF\ : \ \ \Om \by \Big(\R^N\by \R^{Nn}\by \R^{Nn^2}_s \by \cdots \by \R^{Nn^p}_s \Big) \larrow \R^M
\eeq
a Carath\'eodory map. The theory we propose herein applies to measurable solutions $u:\R^n \supseteq \Om \larrow \R^N$ of the system
\beq \label{1.2}
\mF\Big(x,u(x),\D u(x),...,\D^pu(x)\Big)\hspace{1pt} =\hspace{1pt} 0, \quad x\in\Om,
\eeq
without any further restrictions on $\mF$ and $u$. In \eqref{1.1}-\eqref{1.2}, $\R^{Nn}$ symbolises the space of $N\!\by \!n$ matrices and $\R^{Nn^p}_s$ symbolises the space of symmetric tensors 
\[
\begin{split}
&\Big\{ \X \in \R^{Nn^p}\ \big| \ \X_{\al i_1 ... i_p} = \X_{\al \si(i_1 ... i_p)}\ , \ \al =1,...,N, \\
 & \ i_k =1,...,n,\ k=1,...,p,\ \si  \text{ permutation of size }p \Big\}
 \end{split}
\]
wherein the gradient matrix $\D u = ( \D_i u_\al(x))_{i=1,...,n}^{\al=1,...,N}$ and the $p$th order derivative 
\[
\D^p u(x)\ =\ \left( \D^p_{i_1 ... i_p}u_\al(x)\right)_{i_1,...,i_p \in \{1,...,n\}}^{\al=1,...,N}
\]
of (smooth) maps are respectively valued. Evidently, $\D_i= \p/\p x_i$, $x=(x_1,...,x_n)^\top$, $u=(u_1,...,u_N)^\top$ and $\R^{Nn^1}_s=\R^{Nn}$. Since we will not assume that the solutions are locally integrable on $\Om$, the derivatives $\D u$, ..., $\D^p u$ may not have classical meaning, \emph{not even in the sense of distributions}.

The starting point of our approach in not based either on duality or on the maximum principle. Instead, it builds on the probabilistic representation of infinitesimal limits of difference quotients by using \emph{Young measures}. This concept was introduced in the 1930s (\cite{Y}) in order to show existence of ``relaxed" solutions to nonconvex variational problems for which the minimum may not be attained. Nowadays Young measures form a full-blown active area of general topology (\cite{CFV, FG, V}), whilst their utility in Calculus of Variations and PDE theory renders them  indispensable tools for applications (\cite{E, P, FL, M, DPM, KR}), especially in the quantification of the failure of strong convergence due to oscillations and/or concentrations. 

In the present framework, Young measures valued into \emph{compact tori and spheres} (instead of Euclidean spaces as in the aforementioned applications) are utilised in order to define generalised solutions of \eqref{1.2} by applying them to the \emph{difference quotients} of the candidate solution. The notion is pedagogically derived later, but the idea of the definition when $p=1$ in \eqref{1.1} can be briefly motivated as follows: let $u : \R^n \supseteq \Om \larrow \R^N$ be a strong $W^{1,\infty}(\Om,\R^N)$ solution to
\beq  \label{1.3}
\mF\big(\cdot,u,\D u\big)\hspace{1pt}=\hspace{1pt} 0, \quad \text{a.e.\ on }\Om.
\eeq
We aim at finding a ``weak" formulation of \eqref{1.3} which makes sense when $u$ is merely measurable. To this end, we restate \eqref{1.3} as
\beq \label{1.4}
\int_{\R^{Nn}} \Phi(X)\hspace{1pt} \mF\big( x,u(x),X\big)\hspace{1pt} d[\de_{\D u(x)} ](X)\hspace{1pt} =\hspace{1pt} 0, \quad \text{ a.e. }x\in \Om,
\eeq
for any compactly supported $\Phi \in C_c\big( \R^{Nn} \big)$. Namely, we switch from the classical viewpoint of the gradient as a map $\D u : \Om\sub\R^n \larrow \R^{Nn}$ by seeing it a probability-valued map given by the Dirac mass at $\D u$:
\[
\ \ \de_{\D u}\ : \ \ \Om\sub\R^n \larrow \mathscr{P}(\R^{Nn}),\quad    x\lmapsto \de_{\D u(x)}.
\]
Further, we may restate that $\D u$ is the limit in measure of the difference quotients $\D^{1,h}u$ as $h\ri 0$ by writing
\beq \label{1.5}
\de_{\D^{1,h}u} \weakstar\hspace{1pt} \de_{\D u}, \quad \text{as }h\ri0.
\eeq
The weak* convergence above is meant in the Young measures valued into $\R^{Nn}$, that is the set of measurable probability-valued mappings $\Om\sub\R^n \larrow \mathscr{P}(\R^{Nn})$ (for details we refer to Section \ref{section2}). The rationale of the reformulation \eqref{1.4}-\eqref{1.5} of \eqref{1.3} is that \emph{we may thus allow for general probability-valued ``diffuse gradients" of measurable maps which may not be concentration measures.} This is indeed possible if we replace $\R^{Nn}$ by its $1$-point spherical compactification $\smash{\overline{\R}}^{Nn} := \R^{Nn}\cup\{\infty\}$. 
By considering instead the maps $(\de_{\D^{1,h}u})_{h\neq 0}$ as \emph{Young measures valued into $\smash{\overline{\R}}^{Nn}$}, we obtain the necessary compactness and  
we \emph{always} have subsequential weak* limits in the set of Young measures $\Om\sub\R^n \larrow \mathscr{P}\big(\smash{\overline{\R}}^{Nn}\big)$:
\beq  \label{1.5a}
\de_{\D^{1,h_i}u} \weakstar\hspace{1pt} \mD u, \quad \text{as }h_i\ri0.
\eeq
Then, we interpret \eqref{1.3} for just measurable maps $u : \R^n \supseteq \Om \larrow \R^N$ as
\beq  \label{1.5b}
\int_{\smash{\overline{\R}}^{Nn}} \Phi(X)\hspace{1pt} \mF\big( x,u(x),X\big)\hspace{1pt} d[\mD u(x)](X)\hspace{1pt} =\hspace{1pt} 0, \quad \text{ a.e. }x\in \Om,
\eeq
 for any ``test function" $\Phi \in C_c\big( \R^{Nn} \big)$ and any ``diffuse gradient" $\mD u$. Up to a minor technical adaptation (we may need to expand derivatives with respect to non-standard frames determined by $\mF$) \eqref{1.5a} and \eqref{1.5b} essentially constitute the definition of \textbf{diffuse derivatives} and \textbf{$\mD$-solutions} in the special case of \eqref{1.3} and will be the central notion of solution in this paper.

Our motivation to introduce and study generalised solutions for nonlinear PDE systems is primarily sparked by the recently discovered systems associated to vectorial Calculus of Variations in the space $L^\infty$ and in particular the model $\infty$-Laplace system. Calculus of Variations in $L^\infty$ has a long history which started in the 1960s (\cite{A1}-\cite{A5}). Aronsson was the first to consider variational problems for the supremal functional 
\beq  \label{1.6}
\ \ \ \ \ E_\infty(u,\Om')\hspace{1pt} :=\hspace{1pt} \big\| H(\cdot,u,\D u) \big\|_{L^\infty(\Om')}, \quad u\in W^{1,\infty}(\Om,\R^N), \ \Om'\Subset \Om.
\eeq
He studied the case $N=1$ and introduced the appropriate $L^\infty$-notion of minimisers, derived the PDE which is the $L^\infty$-analogue of the Euler-Lagrange equation and studied its classical solutions. In the simplest case of $H(p)=|p|^2$, the $L^\infty$-equation is called the $\infty$-Laplacian and reads
\beq \label{1.7}
\De_\infty u\hspace{1pt} :=\hspace{1pt} \D u \ot \D u : \D^2 u\hspace{1pt} =\hspace{1pt} 0.
\eeq
Since then, the field has undergone huge development due to both the intrinsic mathematical interest and the important for applications: minimisation of the maximum provides more \emph{realistic} models when compared to the classical case of integral functionals. A basic difficulty is that \eqref{1.7} possesses singular solutions. Aronsson himself exhibited this in \cite{A6,A7} and the field remained dormant until the 1990s when the development of viscosity solutions led to an explosion of interest (e.g.\ \cite{C, BEJ, E, E2} and for a pedagogical introduction see \cite{K8}).

Until recently, the study of supremal functionals in conjunction to their associated PDEs was essentially restricted to $N=1$. The principal obstruction appears to be the absence of an efficient theory of generalised solutions allowing the study of general systems, including those arising in $L^\infty$. For instance, the deep contributions \cite{BJW1, BJW2} essentially aimed at studying only the functional when $N\geq 2$. The foundations of the vector case, including the discovery of the appropriate system counterpart of \eqref{1.7}, the correct vectorial $L^\infty$-minimality notion and the study of classical solutions have been laid in a series of recent papers of the author (\cite{K1}-\cite{K6}). In the model case of 
\beq  \label{1.8}
\ \ \ \ \ E_\infty(u,\Om')\hspace{1pt} =\hspace{1pt} \big\| |\D u|^2 \big\|_{L^\infty(\Om')}, \quad u\in W^{1,\infty}(\Om,\R^N), \ \Om'\Subset \Om
\eeq
(where $|\D u|$ is the Euclidean norm of the gradient on $\R^{Nn}$), the analogue of the Euler-Lagrange equation is the $\infty$-Laplace system:
\beq \label{1.9}
\De_\infty u \hspace{1pt} :=\hspace{1pt} \Big(\D u \ot \D u + |\D u|^2[\![ \D u ]\!]^\bot \! \ot \mathrm{I} \Big):\D^2 u\hspace{1pt} =\hspace{1pt} 0.
\eeq
In the above, $[\![ \D u(x)]\!]^\bot$ denotes the orthogonal projection on the orthogonal complement of the range of the matrix $\D u(x)$. In index form \eqref{1.9} reads
\[
\ \ \ \sum_{\be=1}^N\sum_{i,j=1}^n \Big(\D_i u_\al \hspace{1pt} \D_ju_\be +\hspace{1pt} |\D u|^2 [\![\D u]\!]_{\al \be}^\bot \hspace{1pt} \de_{ij}\Big)\hspace{1pt} \D^2_{ij}u_\be\hspace{1pt} =\hspace{1pt} 0, \ \ \ \ \al=1,...,N 
\]
and $[\![ \D u ]\!]^\bot = \hspace{1pt} \text{Proj}_{(R(\D u))^\bot}$. An additional difficulty of \eqref{1.9} which is not present in the scalar case of \eqref{1.7} is that the nonlinear operator may have \emph{discontinuous coefficients} even when applied to smooth maps because the new term involving $[\![ \D u(x)]\!]^\bot$ depends on the dimension of the tangent space of $u(\Om)$ at $x$ (\cite{K1,K6}). Let us also note that almost simultaneously to \cite{K1}, Sheffield and Smart \cite{SS} studied the relevant problem of vectorial optimal Lipschitz extensions and derived a different singular version of ``$\infty$-Laplacian", which in the present setting amounts to changing in \eqref{1.8} from the Euclidean to the operator norm on $\R^{Nn}$.

A further motivation to introduce generalised solutions stems from the insufficiency of the current PDE approaches to handle even elliptic linear systems with rough coefficients. For example, if $\A$ is a continuous symmetric 4th order tensor on $\R^{Nn}$ satisfying the strict Legendre-Hadamard condition, for the divergence system
\[
\sum_{\be=1}^N\sum_{i,j=1}^n \D_i \Big(\A_{\al i \be j}(x)\hspace{1pt} \D_j u_\be(x)\Big)\hspace{1pt} =\hspace{1pt} 0,  \quad \al=1,...,N,
\]
``everything" is known: existence-uniqueness of weak solutions, regularity, etc (see e.g. \cite{GM}). On the other hand, for its non-divergence counterpart
\beq \label{1.12AA}
\sum_{\be=1}^N\sum_{i,j=1}^n \A_{\al i \be j}(x)\hspace{1pt} \D^2_{ij}u_\be(x)\hspace{1pt} =\hspace{1pt} 0,  \quad \al=1,...,N,
\eeq
``nothing" is known, not even what is a meaningful notion of generalised solution, unless $\A$ is $C^{0,\al}$ and \emph{strictly elliptic} in which case a priori estimates guarantee that solutions of \eqref{1.12AA}, if they exist, must be classical (\cite{GM}). To the best of our knowledge there are \emph{no results} for \eqref{1.12AA} in the general case. If $\A$ is monotone (i.e.\ $\A_{\al i \be j}=\de_{\al \be}A_{ij}$), the system decouples to $N$ independent equations and can be treated in the viscosity sense. 

In the present paper, after motivating, introducing and developing some basic theory of $\mD$-solutions for general systems (Section \ref{section2}), we apply it to two important problems. Accordingly, we first consider the Dirichlet problem 
\beq \label{1.10}
\left\{
\begin{array}{rl}
\De_\infty u \hspace{1pt} =\hspace{1pt} 0, & \text{ in }\Om,\\
u\hspace{1pt} =\hspace{1pt} g,  & \text{ on }\p\Om,
 \end{array}
 \right.
\eeq
when $\Om \sub \R^n$ is an open domain with finite measure, $n=N$ and $g\in W^{1,\infty}(\Om,\R^n)$. In Section \ref{section3} we prove existence of $\mD$-solutions to \eqref{1.10} in $W_g^{1,\infty}(\Om,\R^n)$ with extra geometric properties (Theorem \ref{theorem20}, Corollary \ref{corollary21}). \emph{The question of uniqueness for vectorial $L^\infty$ problems has already been answered negatively in \cite{K2}} even in the class of \textbf{smooth} solutions (Remark \ref{rem30}). 

The idea of the proof has two main steps (see Subsection \ref{subsection3.1}). We first apply the Dacorogna-Marcellini Baire Category method (\cite{DM}) which is an analytic counterpart of Gromov's Convex Integration and prove existence of a $W^{1,\infty}$ map solving a first order differential inclusion associated to \eqref{1.10}. Next, we characterise this map as a $\mD$-solution to \eqref{1.10} by utilising the machinery of Section \ref{section2}. Along the way we establish a general tool of independent interest which goes far beyond the $\infty$-Laplacian and provides a method of constructing $\mD$-solutions to ``tangent equations" (Theorem \ref{lemma31}).

The second main question we consider in this paper concerns the existence, uniqueness and (partial) regularity of $\mD$-solutions to the Dirichlet problem
\beq \label{1.11}
\left\{
\begin{array}{rl}
\mF(\cdot,\D^2 u) \hspace{1pt}=\hspace{1pt} f, & \text{ in }\Om,\\
u\hspace{1pt}=\hspace{1pt} 0, & \text{ on }\p\Om,
 \end{array}
 \right.
\eeq
when $\Om\Subset \R^n$ is a $C^2$ convex domain,  $F : \Om \by \R^{Nn^2}_s\larrow \R^N$ is a Carath\'eodory map and $f\in L^2(\Om, \R^N)$. The essential hypothesis guaranteeing well posedness is a \emph{degenerate ellipticity} condition which requires $\mF$ to be ``controllably away" from a degenerate linear operator (Definition \ref{K-Condition}). Our condition is relatively strong, but classical examples (see e.g.\ \cite{LU}) show that \emph{even in the scalar case}, the Dirichlet problem for the uniformly elliptic equation $A(x):\D^2 u(x)=f(x)$ is \emph{not well posed} if $A$ is discontinuous and extra conditions are required. \eqref{1.11} has first been considered by Campanato \cite{C1}-\cite{C3} under a strict ellipticity assumption of Cordes type which implies \eqref{1.11} is well posed in $(W^{2,2}\cap W^{1,2}_0)(\Om,\R^N)$. Very recently, the author (\cite{K9, K11} and \cite{K7}) generalised these results by proving well posedness in the same space under a weaker condition. The latter results for strong solutions of strictly elliptic systems were stepping stones to the approach we develop herein for $\mD$-solutions of degenerate systems. 

In Section \ref{section4} we prove existence of a unique $\mD$-solution to \eqref{1.11}  which satisfies a \emph{new type of partial regularity}, possessing differentiable projections only along certain rank-one lines (Theorem \ref{theorem30}). This regularity is \emph{optimal} (Remark \ref{all remarks}). In particular, the solution may not be even $W^{1,1}_{\text{loc}}$ and does not enjoy any conventional partial regularity of the type of being more regular on a set of full measure. An extra difficulty is the satisfaction of the boundary condition since under this low regularity there is no trace operator. 

The proof is rather long and is based on the study of \eqref{1.11} for linear degenerate systems with constant coefficients  \emph{in the $\mD$-sense} and on a ``perturbation device". The solvability of the linear problem involves approximation and a priori \emph{partial estimates} (Theorem \ref{theorem31}). Well posedness of \eqref{1.11} is established via fixed point in an appropriate functional ``fibre space"  tailored to the degenerate case (\eqref{4.4A},\eqref{4.5A}). The fibre space is an extension of the classical Sobolev space and consists of \emph{partially} regular maps being weakly differentiable only along certain ``elliptic" rank-one directions. We then characterise the fixed point in the fibre space as the unique $\mD$-solution of \eqref{1.11}.   
  
We conclude this introduction by noting that in the companion paper \cite{K12} of this work we have utilised the present framework to prove existence and partial regularity of absolutely minimising $\mD$-solutions to the system of equations arising from \eqref{1.6} when $n=1$ (see also the joint works \cite{AK, KP1} with Abugirda and Pryer). Moreover, in the most recent papers \cite{K13}-\cite{K15} and also jointly with Croce, Pisante and Pryer in \cite{CKP, KP2} we have obtained various explorative results by utilising $\mD$-solutions. We hope that the systematic theory proposed herein will be the starting point for future developments.

\section{Theory of $\mD$-solutions for fully nonlinear systems} \label{section2}

\subsection{Preliminaries}  \label{subsection2.1} We begin with some introductory material which will be used throughout freely, perhaps without explicit reference to this subsection. 
\ms

\noi \textbf{Basics.} Let $n,N\in \N$ be fixed, which in this paper will always be the dimensions of domain and range respectively of our candidate solutions $u:\R^n \supseteq \Om \larrow \R^N$. By $\Om$ we will always mean an open subset of $\R^n$. Unless indicated otherwise, Greek indices $\al,\be,\ga,...$ will run in $\{1,...,N\}$ and latin indices $i,j,k,...$ (perhaps indexed $i_1,i_2,...$) will run in $\{1,...,n\}$, even when the range is not given explicitly. The norms $|\cdot|$ appearing throughout will always be the Euclidean, while the Euclidean inner products will be denoted by either ``$\cdot$" on $\R^n,\R^N$ or by ``$:$" on tensor spaces, e.g.\ on $\R^{Nn}$ we have $ |X|^2 = \sum_{\al,i}X_{\al i}X_{\al i}$ $\equiv X\!:\!X$ and on $\R^{Nn^2}_s$ we have $|\X|^2 = \sum_{\al,i,j}\X_{\al i j}\X_{\al i j}$  $\equiv \X \!:\!\X$,
etc. The standard bases on $\R^n$, $\R^N$, $\R^{Nn}$ will be denoted by $\{e^i\}$, $\{e^\al \}$ and $\{ e^\al \ot e^i\}$. By introducing the \emph{symmetrised tensor product}
\beq \label{2.1}
a \vee b \hspace{1pt} :=\hspace{1pt} \frac{1}{2}\Big(a\ot b \hspace{1pt} +\hspace{1pt} b\ot a \Big), \ \ \ a,b\in \R^n,
\eeq
we will write $\big\{ e^\al \ot ( e^{i_1} \vee ... \vee e^{i_p} ) \big\}$ for the standard basis of the $\R^{Nn^p}_s$. We will follow the convention of denoting vector subspaces of Euclidean spaces as well as the orthogonal projections on them by the same symbol. For example, if $\Si \sub \R^N$ is a subspace, we denote the projection map $\text{Proj}_\Si : \R^N \larrow \R^N$ by just $\Si$ and we have $\Si^2=\Si^\top =\Si \in \R^{N^2}_s$.  We will also systematically use the Alexandroff $1$-point compactification $\R^{Nn^p}_s \cup \{\infty\}$ of the space $\R^{Nn^p}_s$. Its metric will be the standard one which makes it homeomorphic to the sphere of the same dimension (via the stereographic projection which identifies $\{\infty\}$ with the north pole). We will denote it by $\smash{\overline{\R}}^{Nn^p}_s$. We note that all balls and distances taken in $\R^{Nn^p}_s$ (which we will view as a metric vector space isometrically contained into $\smash{\overline{\R}}^{Nn^p}_s$) will be the Euclidean. Similar consideration apply to the torus $\smash{\overline{\R}}^{Nn} \by ... \by \smash{\overline{\R}}^{Nn^p}_s$ and its densely and compactly contained metric vector space $\R^{Nn} \by ... \by  \R^{Nn^p}_s$. Our measure theoretic and function space notation is either standard as e.g.\ in \cite{E2,D, EG} or self-explanatory. For example, the modifier ``measurable" will always mean ``Lebesgue measurable", the Lebesgue measure on $\R^n$ will be denoted by $|\cdot|$, the $s$-Hausdorff measure by $\mH^s$, the characteristic function of a set $A$ by $\chi_A$, the standard $L^p$ and Sobolev spaces of maps  $u : \R^n \supseteq \Om \larrow \Si \sub \R^N$ by $L^p(\Om,\Si), W^{m,p}(\Om,\Si)$, etc. 

\ms

\noi \textbf{General frames, derivative expansions, difference quotients.} In what follows we will need to consider non-standard orthonormal frames of $\R^{Nn^p}_s$ and express derivatives $\D^p u$ with respect to them. Let $\{E^1,...,E^N\}$ be an orthonormal frame of $\R^N$ and suppose that for each $\al=1,...,N$ we have an orthonormal frame $\{E^{(\al)1},...,E^{(\al)n}\}$ of $\R^n$. For these orthonormal frames, we will  equip the spaces $\R^{Nn}=\spn[ \big\{ E^{\al i}  \big\}]$ and $\R^{Nn^p}_s=\spn[ \big\{ E^{\al i_1... i_p}  \big\}]$ with: 
\beq \label{2.2}
E^{\al i} := \hspace{1pt} E^\al \ot E^{(\al)i}, \ \ \ \ \ E^{\al i_1... i_p}:= \hspace{1pt}E^\al \ot \left( E^{(\al)i_1} \vee ... \vee E^{(\al)i_p} \right).
\eeq
For these frames, let $\D_{E^{(\al)i}}$ and $\D^p_{ E^{(\al)i_p} ... E^{(\al)i_1}} \!=\D_{E^{(\al)i_p}}\cdots \D_{E^{(\al)i_1}}$ denote the directional derivatives of $1$st and $p$th order along the respective directions. The gradient $\D u$ of a map $u:\R^n \supseteq \Om \larrow \R^N$ can be expressed as 
\beq \label{2.3a}
\D u\hspace{1pt} =\hspace{1pt} \sum_{\al, i} \Big( E^{\al i} : \D u \Big) \hspace{1pt} E^{\al i}\hspace{1pt} =\hspace{1pt} \sum_{\al, i}\Big( \D_{E^{(\al)i}}(E^\al \cdot u)\Big) \hspace{1pt} E^{\al i}
\eeq
and in general the $p$th order derivative $\D^p u$ as
\beq  \label{2.3b}
\begin{split}
\D^pu\hspace{1pt} &=\!\! \sum_{\al, i_1,...,i_p} \Big( E^{\al i_1...i_p} : \D^p u \Big) \hspace{1pt} E^{\al i_1...i_p}  = \!\! \sum_{\al, i_1,...,i_p} \!\! \Big( \D^p_{E^{(\al)i_1} ...E^{(\al)i_p}}(E^\al \cdot u)\Big) \hspace{1pt} E^{\al i_1...i_p} .
\end{split}
\eeq
We will also use the following notation for the \emph{$p$th order Jet of $u$}:
\[
\D^{[p]} u\hspace{1pt} :=\hspace{1pt} \big(\D u,\D^2 u,...,\D^pu\big).
\]
Given $a\in \R^n$ with $|a|=1$ and $h\in \R\set \{0\}$, when $x,x+ah \in \Om$ the $1$st difference quotient of $u$ along the direction $a$ at $x$ will be denoted by
\beq  \label{2.4}
\D^{1,h}_a u(x)\hspace{1pt} :=\hspace{1pt} \frac{u(x+ha)-u(x)}{h}.
\eeq 
By iteration, if $h_1,...,h_p\neq 0$ the $p$th order difference quotient along $a_1,...,a_p$ is
\beq  \label{2.5}
\D^{p,h_p...h_1}_{a_p...a_1} u\hspace{1pt} :=\hspace{1pt} \D^{1,h_p}_{a_p}\Big( \cdots \big( \D^{1,h_1}_{a_1}u\big) \Big).
\eeq

\ms

\noi \textbf{Young Measures into compact spaces.} This subsection collects basic material that can be found in different guises and greater generality e.g.\ in \cite{CFV, FG, V}. Let $E\sub \R^n$ be a measurable set and $\mK$ a compact subset of some Euclidean space $\R^d$. Later we will take $\mK$ to be either the sphere $\smash{\overline{\R}}^{Nn^p}_s$ or the torus $\smash{\overline{\R}}^{Nn}\! \by \cdots \by \smash{\overline{\R}}^{Nn^p}_s$. Consider the space $L^1\big( E, C(\mK)\big)$ of strongly measurable maps in the standard Bochner sense, where $C(\mK)$ is the space of real continuous functions on $\mK$ (for details we refer e.g.\ to \cite{Ed, FL, F} and references therein). The elements of $L^1\big( E, C(\mK)\big)$ are Carath\'eodory functions $\Phi : E \by \mK \larrow \R$ (i.e.\ $\Phi(\cdot,X)$ is measurable for $X\in \mK$ and $\Phi(\cdot,X)$ is continuous for a.e.\ $x\in E$) for which $\| \Phi \|_{L^1( E, C(\mK))}:= \int_E \max_{X\in \mK} \big|\Phi(x,X)\big|\hspace{1pt} dx < \infty$. It is well-known that (see e.g.\ \cite{FL}) that
\[
\left( L^1\big( E, C(\mK)\big) \right)^* \hspace{1pt} =\hspace{1pt} L^\infty_{w^*}\big( E,\mM(\mK) \big).
\]
The dual space above consists of measure-valued maps $E \ni x \longmapsto \vartheta(x) \in \mM(\mK)$ which are weakly* measurable, that is for any fixed open set $\mathcal{U} \sub \mK$, the function $E\ni x\lmapsto [\vartheta(x)](\mathcal{U}) \in \R$ is measurable. The norm of the space is $\| \vartheta \|_{L^\infty_{w^*} ( E,\mM(\mK) )} := {\ess\hspace{1pt}\sup}_{x\in E} \left\|\vartheta(x) \right\|$, where ``$\|\cdot\|$" denotes the total variation. Since $L^1\big( E, C(\mK)\big)$ is separable, the closed unit ball of $L^\infty_{w^*}\big( E,\mM(\mK) \big)$ is sequentially weakly* compact.  The duality pairing $\langle\cdot,\cdot\rangle :  \ L^\infty_{w^*}\big( E,\mM(\mK) \big) \by L^1\big( E, C(\mK)\big) \larrow \R$ is given by
\[
\langle \vartheta, \Phi \rangle\hspace{1pt} :=\hspace{1pt} \int_E \int_{\mK} \Phi(x,X)\hspace{1pt} d[\vartheta(x)] (X)\hspace{1pt} dx.
\]

\begin{definition}[Young Measures] The subset of the unit sphere of $L^\infty_{w^*}\big( E,\mM(\mK) \big)$ which consists of probability-valued maps is called the set of Young measures:
\[
\mY(E,\mK)\hspace{1pt} :=\hspace{1pt} \Big\{ \vartheta\hspace{1pt} \in \hspace{1pt} L^\infty_{w^*}\big( E,\mM(\mK) \big)\hspace{1pt} : \hspace{1pt} \vartheta(x) \in \mP(\mK),\text{ for a.e. }x\in E\Big\}.
\]
\end{definition}

\begin{remark}[Properties of $\mY(E,\mK)$] \label{remark2} The following well known facts will be extensively used hereafter (for proofs see e.g.\ \cite{FG}): 
\smallskip

\noi (i) [\textbf{weak* compactness}] The set $\mY(E,\mK)$ is convex and (by the compactness of $\mK$, it can be shown) it is sequentially weakly* compact in $L^\infty_{w^*}\big( E,\mM(\mK) \big)$. Hence, for any $(\vartheta^m)_1^\infty$, there is a $\vartheta$ and a subsequence along which $\vartheta^{m_j}\weakstar \vartheta$ as $j\ri \infty$.

\smallskip
\noi (ii)  [\textbf{Young measures induced by functions}] Every measurable map $v : E\sub \R^n \larrow \mK$ induces a Young measure $\de_v \in \mY(E,\mK)$ given by $\de_v(x):= \de_{v(x)}$. 

\smallskip

\noi (iii)  [\textbf{weak* LSC}] We have the following one-sided characterisation of weak* convergence:  $\vartheta^m\weakstar \vartheta$ in $\mY(E,\mK)$ if and only if $\langle \vartheta,\Psi\rangle \hspace{1pt} \leq\hspace{1pt} {\lim \inf}_{m\ri \infty} \langle \vartheta^m,\Psi\rangle$ for any bounded from below function $\Psi : E \by \mK \hspace{1pt} \larrow \hspace{1pt} (-\infty,+\infty]$ measurable in $x$ for all $X\in \mK$ and lower semicontinuous in $X$ for a.e.\ $x\in E$.
\end{remark}

The next result is a minor variant of a classical result which we give together with its short proof because it plays a fundamental role in our setting.

\begin{lemma} \label{lemma2} Suppose $E\sub \R^n $ is measurable and $v^m,v^\infty : E\larrow \mK$ are measurable maps, $m\in \N$. Then, up to the passage to subsequences, we have $v^m \larrow v^\infty$ a.e.\ on $E$ if and only if $\de_{v^{m}} \weakstar \de_{v^\infty}$ in $\mY(E,\mK)$.
\end{lemma}

\BPL \ref{lemma2}.  ($\Rightarrow$) If $v^m \larrow v^\infty$ a.e.\ on $E$, by Remark \ref{remark2} there is $(v^{m_k})_1^\infty$ such that $\de_{v^{m_k}}\weakstar \vartheta^\infty$ in $\mY(E,\mK)$. If $\Phi \in L^1\big( E, C(\mK)\big)$, we have
\[
\int_{E} \Phi \big(x,v^{m_k}(x)\big)\hspace{1pt} dx\hspace{1pt} \larrow \hspace{1pt} \int_{E}\int_\mK \Phi(x,X)\hspace{1pt} d[\vartheta^\infty(x)](X)\hspace{1pt} dx
\]
and also, the $L^1$ bound $|\Phi (\cdot,v^{m_k}) | \leq \max_{X\in \mK} |\Phi(\cdot,X)|$ gives $\Phi (\cdot,v^{m_k}) \larrow \Phi (\cdot,v^{\infty})$ in $L^1(E)$. Hence, by uniqueness of limits $\vartheta^\infty=\de_{v^\infty}$ a.e.\ on $E$.

\smallskip

($\Leftarrow$)  If $\de_{v^{m}} \weakstar \de_{v^\infty}$ in $\mY(E,\mK)$, we choose $\Phi(x,X):= |X-v^\infty(x)|$ where $|\cdot|$ denotes the norm of $\R^d$ restricted to the compact set $\mK$. Then, for any $\e>0$
\[
\begin{split}
0\hspace{1pt}= \int_E \Phi  (\cdot,v^\infty ) \hspace{1pt} = \lim_{m\ri \infty}  \int_E \Phi  (\cdot,v^m )\hspace{1pt} \geq\hspace{1pt} \e\hspace{1pt} \limsup_{m\ri \infty}  \Big|\big\{  |v^m -v^\infty | >\e \big\} \Big|.
\end{split}
\]
Hence, $v^{m} \larrow v^\infty$ in measure on $E$ which gives $v^{m_l} \larrow v^\infty$ a.e.\ on $E$.
 \qed
 
 \ms

\subsection{Motivation of the notions} \label{subsection2.2} We seek to find a meaningful notion of generalised solution for fully nonlinear PDE systems which does not require any more regularity apart from measurability. We derive it in the instructive case of $2$nd order systems. Suppose $\mF$ is as in \eqref{1.1} with $p=2$ and suppose $u : \R^n \supseteq \Om \larrow \R^N$ is a $W^{2,1}_{\text{loc}}(\Om,\R^N)$ strong a.e.\ solution to the system
\beq  \label{2.6}
\mF(\cdot,u,\D u,\D^2 u)\hspace{1pt}=\hspace{1pt} 0, \quad \text{ in }\Om.
\eeq
By the standard equivalence between weak and strong derivatives, the difference quotients converge along subsequence a.e.\ on $\Om$ to the weak derivatives. Hence, 
\[
\mF\Big(\cdot,u,\lim_{m\ri \infty} \D^{1,h_{m} }u ,   \lim_{ m',m''\ri \infty}  \D^{2, h_{m'}h_{m''} }u\Big)\hspace{1pt} =\hspace{1pt} 0, 
\]
a.e.\ on $\Om$. Here $\D^{1,h}$, $\D^{2,kh}$ stand for the usual difference quotient operators whose components with respect to standard basis $\D^{1,h}_{e^i}$, $\D^{2,kh}_{e^i e^j}$ are given by \eqref{2.4}, \eqref{2.5}. Since $\mF$ is a Carath\'eodory map, the limits commute with the nonlinearity:
\beq \label{2.8}
 \lim_{ m,m',m''\ri \infty} \mF\Big(\cdot,u,\D^{1,h_{m} }u, \D^{2, h_{m'}h_{m''} }u\Big)\hspace{1pt} =\hspace{1pt} 0, 
\eeq
a.e.\ on $\Om$. The crucial observation is that \eqref{2.8} \emph{is independent of the weak differentiability of $u$} and makes sense if $u$ is merely measurable. How can we represent these limits and turn them into a handy definition? Going back to \eqref{2.6}, we observe that $u$ is a strong solution of \eqref{2.6} if and only if it satisfies
\[
\int_{ {\R}^{Nn} \by\hspace{1pt}  {\R}^{Nn^2}_s} \Phi(X,\X)\hspace{1pt}\mF(\cdot,u,X,\X)\hspace{1pt} d\big[\de_{(\D u,\D^2 u)} \big](X,\X)\hspace{1pt}=\hspace{1pt} 0, \quad \text{ a.e.\ on }\Om,
\]
for any ``test" function $\Phi \in C_c\big( {\R}^{Nn}\! \by {\R}^{Nn^2}_s \big)$. This gives the idea that we can view the difference quotients as Young measures arising from functions, that is $ \de_{\D^{1,h_m}u}: \Om \larrow  \mP\big( \smash{\overline{\R}}^{Nn} \big)$ and $\de_{\D^{2, h_{m'}h_{m''} }u} : \Om \larrow  \mP\big(  \smash{\overline{\R}}^{Nn^2}_s \big)
$. The reason we compactify the space is to obtain weak* compactness. This compensates the possible loss of mass to $\infty$ since the difference quotients may not converge in any classical sense for just measurable maps. However, when considered in the Young measures valued into spheres they have subsequential weak* limits. It is also more fruitful to take limits \emph{separately} (regardless of order), because the resulting object will be a (fibre) product Young measure valued in the compact torus $\smash{\overline{\R}}^{Nn} \! \by \smash{\overline{\R}}^{Nn^2}_s$:
\beq \label{2.9}
\de_{ \big( \D^{1,h_m}u , \D^{2, h_{m'}h_{m''} }u \big)} \weakstar\hspace{1pt} \mD u \by \mD^2 u \ \ \ \text{ in }\ \mathscr{Y}\Big( \Om, \smash{\overline{\R}}^{Nn} \! \by \smash{\overline{\R}}^{Nn^2}_s \Big),
\eeq
subsequentially as $m,m',m'' \ri \infty$ separately.  Then, \eqref{2.8} is equivalent to
\[
\int_{\smash{\overline{\R}}^{Nn} \by\hspace{1pt} \smash{\overline{\R}}^{Nn^2}_s} \Phi(X,\X)\hspace{1pt} \mF\big( \cdot,u,X,\X \big)\hspace{1pt} d\Big[\de_{\big( \D^{1,h_m}u , \D^{2, h_{m'}h_{m''} }u \big)} \Big](X,\X)\hspace{1pt} \larrow \hspace{1pt} 0, 
\]
subsequentially as $m,m',m'' \ri \infty$, a.e.\ on $\Om$, for any $\Phi\in C_c\big( {\R}^{Nn}\! \by {\R}^{Nn^2}_s \big)$. By using Lemma \ref{lemma16A} that follows, we obtain
\[
\int_{\smash{\overline{\R}}^{Nn} \by\hspace{1pt} \smash{\overline{\R}}^{Nn^2}_s} \Phi(X,\X)\hspace{1pt} \mF\big( \cdot,u,X,\X \big)\hspace{1pt} d\big[\mD u \by \mD^2 u \big](X,\X)\hspace{1pt} =\hspace{1pt} 0, \ \ \text{ a.e.\ on } \Om,
\]
 for any $\Phi$. Note that this statement is \emph{independent} of the regularity of the solution. If $u \in W^{1,1}_{\text{loc}}(\Om,\R^N)$ by Lemma \ref{lemma2} we have that $\mD u= \de_{\D u}$ a.e.\ on $\Om$ and the above statement simplifies to  \eqref{2.11AA}. If further $\D^2 u$ exists weakly on $\Om$, by Lemma \ref{lemma2}  we have $\mD ^2 u=\de_{\D^2 u}$ a.e.\ on $\Om$ thus recovering strong solutions.

\subsection{Main definitions and analytic properties}    

We begin by introducing difference quotients taken with respect to frames as in \eqref{2.2}, \eqref{2.3a}, \eqref{2.3b}. The only difficulty is the complexity in the notation so for pedagogical reasons we give the $1$st order case separately from the general $p$th order case.

\begin{definition}[Difference quotients] \label{definition6} Suppose $\{E^1,...,E^N\}$ is an orthonormal frame of $\R^N$ and for each $\al=1,...,N$ we have an orthonormal frame $\{E^{(\al)1},...,E^{(\al)n}\}$ of $\R^n$ whilst the spaces $\R^{Nn^p}_s$ are equipped with the frames of \eqref{2.2}, $p\in \N$. Let $u:\R^n \supseteq \Om \larrow \R^N$ be a measurable map, extended by zero on $\R^n\set\Om$. Given infinitesimal sequences $ (h_m)_{m\in \N}$ and $(h_{\underline{m}} )_{\underline{m} \in \N^p} \sub \big(\R\set\{0\}\big)^{p}$ such that
\[
\begin{split}
\ \ h_m\ri 0 \text{ as }m\ri \infty\ ,\ \ \ h_{\underline{m}}=( h_{m^1},...,h_{m^p} ), \ \ h_{m^q} \ri 0\ \text{ as } m^q\ri\infty,
\end{split}
\]
we define the \textbf{$1$st and $p$th order difference quotients of $u$} (with respect to the fixed reference frames) arising from $(h_m)_{m\in \N}$ and $(h_{\underline{m}} )_{\underline{m} \in \N^p}$ as
\[
\begin{split}
\ \ \ \ & \D^{1,h_m}u \ :\quad \R^n \supseteq \Om \larrow \R^{Nn},  \ \ \ \ \ m\in \N,\\
& \D^{p,h_{\underline{m}}}u \ : \quad \R^n \supseteq \Om \larrow \R^{Nn^p}_s,  \ \ \ \ \underline{m}=(m^1,...,m^p)\in \N^p,
\end{split}
\]
given respectively by
\[
\begin{split}
 \D^{1,h_m}u  \hspace{1pt} &:=  \hspace{1pt} \sum_{\al, i}  \left[ \D^{1,h_m}_{E^{(\al)i} }(E^\al \cdot u) \right] \hspace{1pt} E^{\al i} ,
\\
 \D^{p,h_{\underline{m}}}u  \hspace{1pt} & := \sum_{\al, i_1,...,i_p}  \Big[ \D^{p,h_{m^p}...h_{m^1} }_{E^{(\al)i_p} ...E^{(\al)i_1}}(E^\al \cdot u)\Big]  \hspace{1pt} E^{\al i_1...i_p} .
\end{split}
\]
In the above, the notation in the brackets is as in \eqref{2.4}, \eqref{2.5}. Further, given an infinitesimal sequence with a trigonal matrix of indices
\[
 (h_{\underline{m}} )_{\underline{m} \in \N^{p^2}} \sub \big(\R\set\{0\}\big)^{p^2}, \ \ \ \underline{m}  = 
\left[
\begin{array}{ccc}
m^1_1 & 0          & 0 \ \hspace{1pt} ...\hspace{1pt} \ 0 \\
m_2^1 & m_2^2 & 0 \ \hspace{1pt} ... \hspace{1pt} \ 0 \\
\vdots  &  & \ddots  \ \hspace{1pt} \vdots \\
m_p^1 & m_p^2 & \hspace{1pt}  ...\hspace{1pt} \ \  m_p^p
\end{array}\right] \!, \ \
h_{m_p^q}\ri0\text{ as }m_p^q\ri\infty,
\]
we will denote its nonzero row elements by $\underline{m}_q\hspace{1pt} :=\hspace{1pt} (m_q^1,...,m_q^q) \in \N^q$, $q=1,...,p$. We define the \textbf{$p$th order Jet $\D^{[p],h_{\underline{m}}}u$ of difference quotients of $u$} (with respect to the fixed reference frames) arising from $ (h_{\underline{m}} )_{\underline{m} \in \N^{p^2}} $ as
\[
\D^{[p],h_{\underline{m}}}u  \hspace{1pt} :=\hspace{1pt} \left(\D^{1,h_{\underline{m}_1}}u,... \hspace{1pt} ,\D^{p,h_{\underline{m}_p}}u\right) \ \  :\  \ \ \R^n \supseteq \Om \larrow \R^{Nn} \by \cdots \by\R^{Nn^p}_s.
\]
\end{definition}

\begin{definition}[Multi-indexed convergence] Let $\underline{m}$ be either a vector of indices in $\N^p$ or a lower trigonal matrix of indices in $\N^{p^2}$. The expression $\underline{m}\larrow \infty$ will symbolise separate successive convergence with respect to each entry in the order: 

\smallskip 

\ \ \ \ \ \ \ $m^1\ri \infty$, ..., $m^p\ri \infty$ \hspace{140pt} \ ($\underline{m}$ vector), 

\ \ \ \ \ \ \ $m_1^1\ri \infty$, $m^1_2\ri \infty$, $m_2^2\ri \infty$, ..., $m_p^{p-1}\ri \infty$, $m_p^p\ri \infty$\hspace{1pt} ($\underline{m}$ matrix).
\end{definition}

\begin{definition}[Diffuse derivatives and Jets]  \label{Diffuse Derivatives}
Suppose we have fixed some reference frames as in Definition \ref{definition6}. For any measurable $u : \R^n \supseteq \Om\larrow\R^N$, we define \textbf{diffuse gradients $\mD u$},  \textbf{diffuse $p$th order derivatives $\mD^p u$} and \textbf{diffuse $p$th order Jets $\mD^{[p]} u$ of $u$} as the subsequential limits of difference quotients arising along infinitesimal sequences in the spaces of Young measures valued in the respective spherical/toric compactifications:
\[
\begin{split}
&\hspace{1pt} \de_{\D^{1,h_{m}}u}\weakstar \mD u, \  \ \ \hspace{1pt} \text{ in }\mY\big(\Om,\smash{\overline{\R}}^{Nn}\big), \ \ \quad \quad \quad \quad \quad \ \hspace{1pt} \text{ as }m\ri \infty,
\\
&\hspace{1pt} \de_{\D^{p,h_{\underline{m}}} u}\weakstar \mD^p u, \  \hspace{1pt} \hspace{1pt} \text{ in }\mY\big(\Om,\smash{\overline{\R}}^{Nn^p}_s\big), \ \quad \quad \quad \quad \quad \ \hspace{1pt} \text{ as }\underline{m}\ri \infty ,\  \underline{m}\in \R^{p},\\
&\de_{\D^{[p],h_{\underline{m}}}u}\weakstar \mD^{[p]} u,   \text{ in }\mY\big(\Om, \smash{\overline{\R}}^{Nn} \by \cdots \by\smash{\overline{\R}}^{Nn^p}_s\big), \ \text{ as }\underline{m}\ri \infty,\ \underline{m}\in \R^{p^2} .
\end{split}
\]
\end{definition} 

\begin{remark} As a consequence of the separate convergence, the $p$th order Jet is always a (fibre) product Young measure: $\mD^{[p]}u\ =\ \mD u \by \cdots \by \mD^p u$ .
\end{remark}

We now record that Remark \ref{remark2}(i) implies the existence of diffuse derivatives.

\begin{lemma}[Existence of diffuse derivatives] Every measurable mapping $u : \R^n \supseteq \Om \larrow\R^N$ possesses diffuse derivatives  of all orders, actually at least one for every choice of infinitesimal sequence.
\end{lemma}

\begin{remark}[Nonexistence of distributional derivatives] Since we do not require our maps to be in $L^1_{\text{loc}}(\Om,\R^N)$, they may not possess distributional derivatives.
\end{remark}

In general diffuse derivatives \emph{may not be unique} for nonsmooth maps. However, they are compatible with weak derivatives:

\begin{lemma}[Compatibility of weak and diffuse derivatives] \label{lemma10} If $u \in W^{1,1}_{\text{loc}}(\Om,\R^N)$, then the diffuse gradient $\mD u$ is unique and
$\de_{\D u} = \mD u$, a.e.\ on $\Om$. More generally, if $q\in \{1,...,p-1\}$ and $u \in W^{q,1}_{\text{loc}}(\Om,\R^N)$, then $\mD^{[q]}u$  is unique and 
\[
\mD^{[p]}u\hspace{1pt} =\hspace{1pt} \de_{(\D u,...,\D^{q}u)}\by \mD^{q+1} \by \cdots \by\mD^p u, \ \ \text{ a.e.\ on } \Om.
\]
\end{lemma}

\BPL \ref{lemma10}. It suffice to establish only the $1$st order case. For any fixed $e\in \R^n$ we have $\D^{1,h}_e u \larrow D_e u$  in $L^{1}_{\text{loc}}(\Om,\R^N)$ as $h\ri 0$. We choose $e:=E^{(\al)i}$ and $h:=h_m$ to get $\D^{1,h_m}_{E^{(\al)i}} (E^\al \cdot u)  \larrow \D_{E^{(\al)i}}(E^\al \cdot u)$ in $L^{1}_{\text{loc}}(\Om)$ as $m\ri \infty$. Thus, by \eqref{2.3a}, \eqref{2.3b} and Definition \ref{definition6} we have $\D^{1,h_m}u \larrow \D u$ a.e.\ on $\Om$ as $m\ri \infty$ along a subsequence. Application of Lemma \ref{lemma2} completes the proof.    \qed

\ms

Next we show that the diffuse gradient is a Dirac mass if and only if $u$ is  ``differentiable in measure", a notion introduced and studied by Ambrosio-Mal\'y in \cite{AM}. This notion arose in the study of the regularity of the flow map of ODEs driven by Sobolev vector fields (see \cite{BL}). 

\begin{definition}[Differentiability in measure, cf.\ \cite{AM}] \label{definition12} Let $u : \R^n \supseteq \Om \larrow \R^N$ be measurable. We say that $u$ is differentiable in measure on $\Om$ with derivative the measurable map $\mL D u : \R^n \supseteq \Om \larrow \R^{Nn}$ if for any $\e>0$ and $E\sub \Om$ with $|E|<\infty$,
\[
\lim_{y\ri 0}\left|\left\{ x\in E\hspace{1pt} :\hspace{1pt} \left|\frac{u(x+y) -u(x)-\mL  \D u(x)\hspace{1pt}  y}{|y|}\right|>\e \right\}\right|\hspace{1pt} =\hspace{1pt}0. 
\]
\end{definition}

In \cite{AM} it is shown that this notion is strictly weaker than the classical notion of approximate differentiability (\cite{EG}).

\begin{lemma}[Gradient in measure vs diffuse gradient] \label{lemma13} Let $u : \R^n \supseteq \Om \larrow \R^N$ be measurable and suppose we have fixed some reference frames as in Definition \ref{definition6}.

\noi (a) If $u$ is differentiable is measure with derivative $\mL \D u$, then the diffuse gradient $\mD u \in \mY(\Om,\smash{\overline{\R}}^{Nn})$ is unique and $\mD u = \de_{\mL \D u}$ a.e.\ on $\Om$.

\smallskip

\noi (b) If there exists a measurable map $U : \R^n \supseteq \Om \larrow \R^{Nn}$ such that for any diffuse gradient $\mD u \in \mY(\Om,\smash{\overline{\R}}^{Nn})$ we have $\mD u = \de_{U}$ a.e.\ on $\Om$, then it follows that $u$ is differentiable in measure and $U=\mL \D u$ a.e.\ on $\Om$.
\end{lemma}

\BPL \ref{lemma13}. (a) By choosing $y:=h E^{(\al)i}$ in Definition \ref{definition12} applied to the projection $E^\al \cdot u$ we get that $\D^{1,h}_{ E^{(\al)i} } (E^\al \cdot u)\larrow  E^{\al i}:(\mL \D u )$ as $h\ri 0 $ locally in measure on $\Om$. Thus, for any $h_m\ri 0$, there is $h_{m_k}\ri 0$ such that the convergence is a.e.\ on $\Om$, whence $\mD u =\de_{\mL \D u}$ by Lemma \ref{lemma2}.

(b) We begin by observing a triviality: for any map $f: \R^n \ri \R^N$ we have $f(y)\ri l$ as $y\ri 0$ if and only if for any $y_m \ri 0$, there is $y_{m_k}\ri 0$ such that $f(y_{m_k})\ri l$ as $k\ri \infty$. We continue by noting that by Lemma \ref{lemma2} and our assumption we have that for any $h_m\ri 0$ there is $h_{m_k}\ri 0$ such that $\D^{1,h_{m_k}}u\larrow U$ a.e.\ on $\Om$,  as $k\ri \infty$. Hence, we obtain that $\D^{1,h}u\larrow U$ as $h\ri 0$ (full limit), a.e.\ on $\Om$. Since a.e.\ convergence implies convergence locally in measure, we deduce that $U=\mL \D u$ a.e.\ on $\Om$, as desired.           \qed

\ms

The next notion of solution will be central in this work. For pedagogical reasons, we give it first for $W^{1,1}_{\text{loc}}$ solutions of $2$nd order systems and then in the general case.

\begin{definition}[Weakly differentiable $\mD$-solutions of $2$nd order PDE systems] \label{definition11} Let $\Om \sub \R^n$ be open, $\mF : \Om \by \big(\R^N\by \R^{Nn}\by \R^{Nn^2}_s \big) \larrow \R^M$ a Carath\'eodory map and $u : \R^n \supseteq \Om \larrow \R^N$ a map in $W^{1,1}_{\text{loc}}(\Om,\R^N)$. Suppose we have fixed some reference frames as in Definition \ref{definition6} and consider the PDE system
\beq \label{2.11}
\mF\big(\cdot,u,\D u,\D^2 u\big)\hspace{1pt} =\hspace{1pt}0, \ \ \text{ in }\Om.
\eeq
We will say that $u$ is a \textbf{$\mD$-solution of \eqref{2.11}} when for any diffuse hessian $\mD ^2 u \in \mY(\Om,\smash{\overline{\R}}^{Nn^2}_s)$ of $u$ (Definition \ref{Diffuse Derivatives}) and any $\Phi \in C_c\big( {\R}^{Nn^2}_s \big)$ we have
\beq \label{2.11AA}
\ \ \ \int_{\smash{\overline{\R}}^{Nn^2}_s} \Phi(\X)\hspace{1pt} \mF\big( \cdot,u,\D u,\X \big)\hspace{1pt} d[\mD^2 u ](\X)\hspace{1pt} =\hspace{1pt} 0, \ \ \text{ a.e.\ on } \Om.
\eeq
 \end{definition} 
 
We note that $\mF$ is not actually necessary to be continuous with respect to $(u,\D u)$ and merely Borel measurable suffices. Now we consider the general case. For brevity we will write $\underline{\X} = (\X_1,...,\X_p)$ for the points of the torus $\smash{\overline{\R}}^{Nn} \by\cdots \by \smash{\overline{\R}}^{Nn^p}_s$.

\begin{definition}[$\mD$-solutions for $p$th order PDE systems] \label{definition13} Let $\Om \sub \R^n$ be open and $\mF$ a Carath\'eodory map as in \eqref{1.1}. Suppose $u : \R^n \supseteq \Om \larrow \R^N$ is measurable and we have fixed some reference frames as in Definition \ref{definition6}. Consider the system
\beq \label{2.11a}
\mF\left(x,u(x),\D^{[p]}u(x)\right)\hspace{1pt} =\hspace{1pt}0, \ \ \ x\in \Om.
\eeq
We will say that $u$ is a \textbf{$\mD$-solution of \eqref{2.11a}} when for any $\mD^{[p]}u \in \mY\big(\Om, \hspace{1pt} \smash{\overline{\R}}^{Nn} \! \by\! \cdots$ $\by\! \smash{\overline{\R}}^{Nn^p}_s\big)$ of $u$ (Definition \ref{Diffuse Derivatives}) and any $\Phi \in C_c\big( {\R}^{Nn} \!\by \cdots \by {\R}^{Nn^p}_s \big)$, we have
\[
\ \ \ \ \int_{\smash{\overline{\R}}^{Nn} \by\cdots \by \smash{\overline{\R}}^{Nn^p}_s} \Phi(\underline{\X})\hspace{1pt} \mF\big( x,u(x),\underline{\X}\big)\hspace{1pt} d\big[\mD^{[p]} u (x)\big](\underline{\X})\hspace{1pt} =\hspace{1pt} 0, \quad \text{ a.e. }x\in\Om.
\]
\end{definition}

The following result asserts the fairly obvious fact that $\mD$-solutions and strong solutions are compatible. 

\begin{proposition}[Compatibility of strong with $\mD$-solutions] \label{proposition14} Let $\mF$ a Carath\'eodory map as in \eqref{1.1} and $u \in W^{p,1}_{\text{loc}}(\Om,\R^N)$ (or merely $p$-times differentiable in measure, Definition \ref{definition12}). Then, $u $ is a $\mD$-solution of \eqref{2.11a} on $\Om$ if and only if $u$ is a strong solution  of \eqref{2.11a} a.e.\ on $\Om$.
\end{proposition}

\BPP \ref{proposition14}. It is an immediate consequence of  Lemma \ref{lemma10} (or Lemma \ref{lemma13}) and the motivation of the notions (Subsection \ref{subsection2.2}).   \qed

\ms

Our next result is a simple yet powerful convergence tool which we give in the generality of Young measures and will play an important role in later sections.

\begin{lemma}[Convergence lemma]\label{lemma16A} Suppose that $u^\infty, (u^\mu)_1^\infty$ are measurable maps $ \R^n \supseteq \Om \larrow \R^N$ satisfying $u^\mu \larrow u^\infty$ a.e.\ on $\Om$. Let $\mathbb{W}$ be a finite dimensional metric vector space, isometrically contained into a compactification $\mK$ of $\mathbb{W}$. Suppose we have Carath\'eodory maps $\mF^\infty$, $\mF^\mu$ : $\Om \by \!\big( \R^N \by \mathbb{W} \big) \larrow \R^M$, $\mu\in \N$ such that for a.e.\ $x\in \Om$, $\mF^\mu (x,\cdot , \cdot)  \larrow \mF^\infty(x , \cdot  , \cdot)$ in $C(\R^N\!\by \mathbb{W})$ as $\mu\ri \infty$. Suppose further we have Young measures $\vartheta^\infty,(\vartheta^\mu)_1^\infty \in \mY\big( \Om,\mK\big)$ such that $\vartheta^\mu\weakstar  \vartheta^\infty$ in $\mY\big( \Om, \mK\big)$ as $\mu\ri \infty$. Then, if for a given $\Phi \in C_c(\mathbb{W})$ we have
\[
\int_{\mK} \Phi(\X) \hspace{1pt} \mF^\mu\big(x,u^\mu(x),\X\big)\hspace{1pt} d[\vartheta^\mu(x)](\X)\hspace{1pt} =\hspace{1pt} 0,  \ \text{ a.e. }x\in \Om,
\]
for all $\mu\in \N$, it follows that the same conclusion holds for $\mu=\infty$ as well.
\end{lemma}

\BPL \ref{lemma16A}. We fix $\Phi \in C_c(\mathbb{W})$ and set
\[
\phi^m(x)\hspace{1pt} :=\hspace{1pt} \Big\|\Phi(\cdot)\Big( \mF^m\big(x,u^m(x),\cdot\big)\hspace{1pt}-\hspace{1pt}   \mF^\infty\big(x,u^\infty(x),\cdot\big)\Big) \Big\|_{C(\mathbb{W})}
\] 
and we claim that $\phi^m(x)\larrow 0$ for a.e.\ $x\in\Om$. To see this, fix $x\in \Om$ such that $u^m(x)\larrow u^\infty(x)$ (the set of such points has full measure in $\Om$). Fix also $U\Subset \R^N$ and $W\Subset \mathbb{W}$ such that $u^m(x),u^\infty(x)\in U$ and $\supp(\Phi)\sub W$ for large $m\in \N$. By our assumptions, $\| \mF^m(x,\cdot , \cdot)-\mF^\infty(x,\cdot , \cdot) \|_{C(U \by W)}\larrow 0$ as $m\ri \infty$. If $\om^\infty_x \in C[0,\infty)$ symbolises the modulus of continuity of $U \ni \xi\lmapsto \mF^\infty(x,\xi,\X) \in \R^M$ which is uniform in $\X\in W$,  we have
\[
\begin{split}
|\phi^m(x)|\hspace{1pt} & \leq \hspace{1pt}   \|\Phi\|_{C(W)}\hspace{1pt} \Bigg(  \Big\| \mF^\infty \big( x,u^m(x),\cdot\big)\hspace{1pt}-\hspace{1pt}  \mF^\infty\big(x,u^\infty(x),\cdot\big) \Big\|_{C(W)} \\
&\hspace{60pt} +\hspace{1pt}   \Big\| \mF^m \big( x,u^m(x),\cdot\big)\hspace{1pt}-\hspace{1pt}  \mF^\infty\big(x,u^m(x),\cdot\big) \Big\|_{C(W)}
 \Bigg)
 \end{split}
\]
 \[
\begin{split}
\leq \hspace{1pt}  \|\Phi\|_{C(\mathbb{W})} \Big(\hspace{1pt}  \om^\infty_x\big(\big|u^m(x)-u^\infty(x)\big| \big)\  +\hspace{1pt} \big\| \mF^m(x,\cdot,\cdot)-\mF^\infty(x,\cdot,\cdot) \big\|_{C(U\by W)}\Big),
\end{split}
\]
giving that $|\phi^m(x)| = o(1),$ as $m\ri \infty$. We now fix $R>0$ and set 
\[
{\Om_R}\hspace{1pt} := \hspace{1pt} \Big\{ x\in \Om \hspace{1pt} : \hspace{1pt} \big\|\Phi(\cdot)\hspace{1pt} \mF^\infty\big(x,u^\infty(x),\cdot\big) \big\|_{C(\mathbb{W})} <\hspace{1pt} R   \Big\} \cap \mB_R(0). 
\]
Since $|{\Om_R}|<\infty$, by the Egoroff theorem we can find measurable sets $\{E_i\}_1^\infty \sub {\Om_R}$ such that $|E_i|\ri 0$ as $i\ri \infty$ and for each $i\in\N$ we have $\phi^m\larrow 0 $ in $L^\infty({\Om_R}\set E_i)$ as $m\ri\infty$. Since $|\Om_R|<\infty$, we have $\phi^m\larrow 0 $ in $L^1({\Om_R}\set E_i)$ as well. Further, the functions $\Psi^m(x,\X):=\left|\Phi(\X)\hspace{1pt} \mF^m\big(x,u^m(x),\X\big)\right|$, $m\in \N\cup\{\infty\}$, are elements of $L^1\big({\Om_R}\set E_i, C(\mK)\big)$ because 
\[
\| \Psi^m-\Psi^\infty\|_{L^1({\Om_R}\set E_i, C(\mK))}
\hspace{1pt} \leq \hspace{1pt} \|\phi^m\|_{L^1({\Om_R}\set E_i)}
\]
and for $m$ large we have
\[
\| \Psi^m\|_{L^1({\Om_R}\set E_i, C(\mK))}\hspace{1pt} \leq\hspace{1pt} 1\hspace{1pt}+\hspace{1pt} \| \Psi^\infty\|_{L^1({\Om_R}\set E_i, C(\mK))}\hspace{1pt} \leq\hspace{1pt} 1+\hspace{1pt} |\Om_R| \hspace{1pt} R.
\]
Hence, $\Psi^m \larrow \Psi^\infty$ in $L^1\big({\Om_R}\set E_i, C(\mK)\big)$ and also by assumption $\vartheta^m \weakstar  \vartheta^\infty$ in $\mY\big({\Om_R}\set E_i,\mK \big)$. By the weak*-strong continuity of the pairing
\[
 L^\infty_{w^*}\Big({\Om_R}\set E_i,\mM ( \mK) \Big) \hspace{1pt} \by\hspace{1pt}  L^1\Big( {\Om_R}\set E_i, C (  \mK)\Big) \larrow \R
\]
we may pass to the limit in our hypothesised identity to obtain
\[
\int_{ \mK } \Phi(\X) \mF^\infty\big(x,u^\infty(x),\X\big)\hspace{1pt} d[  \vartheta^\infty (x)](\X)\hspace{1pt}=\hspace{1pt} 0,
\]
for a.e.\ $x\in {\Om_R}\set E_j$. We conclude by letting $j\ri \infty$ and then taking $R\ri \infty$.          \qed

\ms

The next result is a direct consequence of Lemma \ref{lemma16A} and establishes that $\mD$-solutions are well behaved under weak* convergence.

\begin{corollary}[Convergence of $\mD$-solutions] \label{theorem16} Let $(u^\mu)_1^\infty$ be a sequence of maps where each $u^\mu:\R^n \supseteq \Om \larrow \R^N$ is measurable and $u^\mu \larrow u^\infty$ a.e.\ on $\Om$. Let also $(\mF^\mu)_1^\infty$ be  Carath\'eodory maps as in \eqref{1.1}. Suppose each $u^\mu$ is a $\mD$-solution of
\[
\mF^\mu\left(x,u^\mu(x),\D^{[p]}u^\mu(x)\right)\hspace{1pt}=\hspace{1pt} 0, \ \ x\in \Om,
\]
and $  \mF^\mu(x,\cdot,\cdot)    \larrow \mF^\infty(x,\cdot,\cdot)$ uniformly on compact subsets as $\mu\ri\infty$, for a.e.\ $x\in\Om$. If every jet $\mD^{[p]}u^\infty$ can be weakly* approximated by a subsequence of the respective Jets $\mD^{[p]}u^{\mu_\nu}$, then $u^\infty$ is a $\mD$-solution of the limit system for $\mu=\infty$.
\end{corollary}

\begin{remark} \label{remark30} Note that Corollary \ref{theorem16} is \emph{not} a stability result, in the sense that we do not have compactness of diffuse jets as part of the conclusion. In fact, such a result is not possible without extra assumptions which would entail some sort of \emph{a priori estimates}: for instance, consider the sequence $u^\mu(x):=\mu^{-1}\sin(\mu x)$, $x\in\R$. Then, $u^\mu\weakstar u^\infty$ in $W^{1,\infty}(\R)$ where $u^\infty\equiv 0$. However, $\mD u^\mu = \de_{\D u^\mu}\weakstar \vartheta$ in $\mY(\R,\R)$  as $\mu\ri\infty$,  where for a.e.\ $x\in \R$ $\supp(\vartheta(x))=[-1,1]$  
while $\vartheta(x) \neq \mD u^\infty(x) = \de_{\{0\}}$.
\end{remark}

\noi The next result gives equivalent formulations of the definition of $\mD$-solutions. To this end we need some further terminology.  

\begin{definition}[Reduced support] \label{definition24a} Given a probability $\vartheta  \in \mP \big(\smash{\overline{\R}}^{Nn}\by \cdots\by \smash{\overline{\R}}^{Nn^p}_s\big)$, we define its \textbf{reduced support} as
\[
\supp_*(\vartheta)\hspace{1pt} :=\hspace{1pt} \supp(\vartheta) \cap \left(\R^{Nn}\by \cdots \by \R^{Nn^p}_s\right).
\]
\end{definition}

\begin{definition}[Cut offs associated to a map] \label{definition24} Let $u:\R^n \supseteq \Om\larrow \R^N$ be measurable and $\mF$ as in \eqref{1.1}. For any measurable $U : \R^n \supseteq \Om \larrow  {\R}^{Nn}\by \cdots \by{\R}^{Nn^p}_s$ and $R>0$, we define the \textbf{cut off of $U$ associated to $\mF$} as:
\[
[ U]^R\hspace{1pt} :=\hspace{1pt} \left\{
\begin{array}{ll}
U, &\text{ on } \big\{|U| \leq R\big\},\smallskip\\
\textbf{0}^R, & \text{ on } \big\{ |U| > R\big\}.
\end{array}
\right.
\]
Here $\textbf{0}^R$ is a measurable selection of the set-valued mapping 
\[
\Om\ni \hspace{1pt}x\hspace{1pt} \lmapsto\hspace{1pt} \Big\{\mF\big(x,u(x),\cdot\big) =0\Big\}\cap \mB_R(0) \ \sub \hspace{1pt} \Big(\R^{Nn} \by \cdots \by\R^{Nn^p}_s\Big) \set \{\emptyset\},
\]
that is, $\textbf{0}^R : \Om \larrow \R^{Nn} \by \cdots \by\R^{Nn^p}_s$ satisfies $\mF\big(x,u(x),\textbf{0}^R(x)\big)\hspace{1pt}=\hspace{1pt}0$ and $|\textbf{0}^R(x)|\leq R$ for a.e.\ $x\in \Om$.
\end{definition}

The existence of selections is a consequence of Aumann's theorem (see e.g.\ \cite{FL}). If $\mF\big(x,u(x),\cdot\big)$ is linear, we may choose $\textbf{0}^R(x)\equiv 0$ with no $(R,x)$-dependence.

\begin{proposition}[Equivalent definitions for $\mD$-solutions] 
\label{proposition15}  Let $\mF$ be as in \eqref{1.1} and $u:\R^n \supseteq \Om\larrow \R^N$ a measurable map. Then, the following are equivalent:
\begin{enumerate}

\item The map $u$ is a $\mD$-solution of the PDE system
\[
\mF\left(x,u(x),\D^{[p]}u(x) \right)=0, \ \ \ x\in \Om.
\]

\item For any diffuse $p$th order Jet of $u$, we have
\[
\sup_{\underline{\X} \in  \supp_*(\mD^{[p]}u(x)) }\big| \mF\big(x,u(x),\underline{\X} \big)\big|\hspace{1pt} =\hspace{1pt} 0, \ \ \text{ a.e. }x\in \Om.
\]

\item For any diffuse $p$th order Jet of $u$, we have the inclusion
\[
\supp_*\big(\mD^{[p]}u(x)\big) \sub\hspace{1pt} \Big\{\mF\big(x,u(x),\cdot\big)=0 \Big\},\ \ \text{ a.e. }x\in \Om.
\]

\item For any diffuse $p$th order Jet of $u$, we have
\[
\quad \quad \int_{ \R^{Nn} \by \cdots \by\R^{Nn^p}_s} \big| \mF\big(x,u(x),\underline{\X} \big)\big|\hspace{1pt} d\big[\mD^{[p]} u(x) \big](\underline{\X})\hspace{1pt} =\hspace{1pt} 0, \ \ \text{ a.e. }x\in \Om.
\]

\item For any $p$th order Jet of difference quotients of $u$ and any $R>0$, we have
\[
\quad \quad \mF\left(\cdot,u,\big[ \D^{[p],h_{\underline{m}}}u\big]^R\right)\hspace{1pt} \larrow\hspace{1pt} 0,
\]
for a.e.\ $x\in \Om$, as $\underline{m}\ri \infty$ along subsequences.

\item For any $p$th order Jet of difference quotients of $u$ and any $R>0$, we have
\[
\quad \quad \ \  \dist  \Bigg( \big[ \D^{[p],h_{\underline{m}}}u\big]^R(x)\ , \ \Big\{\mF\big(x,u(x),\cdot\big)=0 \Big\} \cap \mB_R(0)\Bigg) \hspace{1pt} \larrow\hspace{1pt} 0,
\]
for a.e.\ $x\in \Om$, as $\underline{m}\ri \infty$ along subsequences.
\end{enumerate}
\end{proposition}

The presence of reduced supports and cut offs is informally interpreted as that the mass which does not escape to infinity actually lies on the zero level set of the coefficients. The proof of Proposition \ref{proposition15} does not rely on the particular structure of diffuse Jets and is a consequence of the next more general result.

\begin{lemma} \label{lemma15a} All the equivalences of Proposition \ref{proposition15} remains true if more generally one replaces the jet $\D^{[p],h_{\underline{m}}}u$ by any measurable sequence
\[
U^m \ :\ \ \R^n \supseteq \Om \larrow \hspace{1pt} {\R}^{Nn}\by \cdots \by{\R}^{Nn^p}_s, \quad m\in \N,
\]
and the respective Jet $\mD^{[p]} u$ by any Young measure $\vartheta \in \mY\left(\Om, \hspace{1pt} \smash{\overline{\R}}^{Nn} \by \cdots \by \smash{\overline{\R}}^{Nn^p}_s\right)$ such that $\de_{U^m}\weakstar \vartheta$ as $m\ri \infty$.
\end{lemma}

\BPL \ref{lemma15a} \& \textbf{Proposition} \ref{proposition15}. We begin by showing (1)$\Leftrightarrow$(3)$\Leftrightarrow$(2) and then we establish that (6)$\Rightarrow$(5)$\Rightarrow$(4)$\Rightarrow$(3)$\Rightarrow$(6).

\ms
\noi (1)$\Rightarrow$(3): Suppose $\de_{U^m}\weakstar \vartheta$ and that for any $\Phi \in C_c\big( {\R}^{Nn}\by \cdots \by{\R}^{Nn^p}_s \big)$ we have 
 \[
\quad \quad \int_{ \smash{\overline{\R}}^{Nn} \by \cdots \by\smash{\overline{\R}}^{Nn^p}_s}\Phi(\underline{\X})\hspace{1pt} \mF\big(x,u(x),\underline{\X} \big)\hspace{1pt} d[\vartheta (x) ](\underline{\X})\hspace{1pt} =\hspace{1pt} 0, \ \ \text{ a.e. }x\in \Om,
\]
whilst for some of these $x\in \Om$ we have $\supp\big(\vartheta_*(x)\big) \not\sub \big\{\mF\big(x,u(x),\cdot\big)=0 \big\}$. Then, there exists $\underline{\X}_0$ with $\mF\big(x,u(x),\underline{\X}_0\big)\neq 0$ and $[\vartheta(x)] \big(\mB_{R}(\underline{\X}_0) \big)>0$ for $R>0$. By continuity, there exist $c_0,R_0>0$ and $\mu \in \{1,...,M\}$ such that $\left|\mF_\mu\big(x,u(x),\cdot\big)\right|\geq c_0$ on $\mB_{R_0}(\underline{\X}_0)$. By choosing $\Phi$ such that $\chi_{\mB_{R_0 /2}(\underline{\X}_0)} \leq  \Phi \leq \chi_{\mB_{R_0}(\underline{\X}_0)}$, we get 
\[
\begin{split}
0\hspace{1pt}&=\hspace{1pt} \left| \int_{{\R}^{Nn}\by \cdots \by{\R}^{Nn^p}_s} \Phi(\underline{\X})\hspace{1pt} \mF_\mu \big(x,u(x),\underline{\X} \big)\hspace{1pt} d[\vartheta(x) ](\underline{\X}) \right|\\
 &=\hspace{1pt} \int_{\mB_{R_0}(\underline{\X}_0)} \Phi(\underline{\X})\hspace{1pt} \big| \mF_\mu\big(x,u(x),\underline{\X} \big)\big|\hspace{1pt} d[\vartheta(x)](\underline{\X}) \\
 & \geq\hspace{1pt} c_0\hspace{1pt} [\vartheta(x)] \big(\mB_{R_0 /2}(\underline{\X}_0) \big).
\end{split}
\]
The above contradiction establishes that the desired inclusion holds a.e.\ on $\Om$.

\ms
\noi (3)$\Rightarrow$(1): Suppose $\supp_*\big(\vartheta(x)\big) \sub \big\{\mF\big(x,u(x),\cdot\big)=0 \big\}$ for a.e.\ $x\in \Om$. Then, for any $\Phi \in C_c\big({\R}^{Nn}\by \cdots \by{\R}^{Nn^p}_s \big)$ and any such $x$, $\Phi(\cdot) \mF\big(x,u(x),\cdot\big)$ vanishes $[\vartheta(x)]$-a.e.\ on $\smash{\overline{\R}}^{Nn}\by \cdots \by\smash{\overline{\R}}^{Nn^p}_s$. Thus, 
\[
\int_{\smash{\overline{\R}}^{Nn}\by \cdots \by\smash{\overline{\R}}^{Nn^p}_s} \Phi(\underline{\X})\hspace{1pt} F \big(x,u(x),\underline{\X} \big)\hspace{1pt} d[\vartheta(x)](\underline{\X})\hspace{1pt} =\hspace{1pt} 0.
\] 

\ms

\noi (3)$\Leftrightarrow$(2): Effectively, they are just restatements of each other and either of them states that for any diffuse $p$th order Jet, for a.e.\ $x\in \Om$ and for all $\underline{\X} \in  \supp_*(\mD^{[p]}u(x))$, we have $\big| \mF\big(x,u(x),\underline{\X} \big)\big|=0$.

\ms

\noi (6)$\Rightarrow$(5): If suffices to show that for a.e.\ $x\in \Om$ and any $R>0$ there is a strictly increasing modulus of continuity $\om_{R,x} \in C[0,\infty)$ such that
\[
\Big|\mF\big(x,u(x),\underline{\X}\big)  \Big|\hspace{1pt} \leq\hspace{1pt} \om_{R,x} \left( \dist\left( \underline{\X}\ , \  \Big\{\mF\big(x,u(x),\cdot\big)=0\Big\} \cap {\mB_{R}(0)} \right)\right),
\]
when $\underline{\X}\in \overline{\mB_{R}(0)}$. In such an event we conclude by choosing $\underline{\X}:=[U^m]^R(x)$. To see the claim, note that for a.e.\ $x\in \Om$ there is such an $\om_{R,x} $ with
\[
\Big|\mF\big(x,u(x),\underline{\X}\big) \hspace{1pt} -\hspace{1pt} \mF\big(x,u(x),\underline{\Y}\big) \Big|\hspace{1pt} \leq\hspace{1pt} \om_{R,x} \big( \big|\underline{\X}-\underline{\Y}\big|\big)
\]
for all $\underline{\X},\underline{\Y}\in \overline{\mB_{R}(0)}$. By choosing $\underline{\Y} \in \big\{\mF\big(x,u(x),\cdot\big) =0\big\}$, we have
\[
\begin{split}
\Big|\mF\big(x,u(x),\underline{\X}\big) \Big|\hspace{1pt} 
&\leq\hspace{1pt} \inf_{\mF(x,u(x),\underline{\Y}) =0, |\underline{\Y}|\leq R} \om_{R,x} \big( \big|\underline{\X}-\underline{\Y}\big|\big)\\
&=\hspace{1pt} \om_{R,x} \left(
\inf_{\mF(x,u(x),\underline{\Y}) =0, |\underline{\Y}|\leq R}  \big|\underline{\X}-\underline{\Y}\big|\right),
\end{split}
\]
as desired.

\ms
\noi (5)$\Rightarrow$(4): We fix $R>0$ and $\Phi \in C_c\big( {\R}^{Nn}\by \cdots \by{\R}^{Nn^p}_s \big)$ such that $\chi_{\mB_{R/2}(0)}\leq \Phi  \leq \chi_{\mB_R(0)}$. For any $k\in \N$, we set
\[
\Om_k\hspace{1pt} :=\hspace{1pt} \left\{ x\in \Om\cap \mB_{k}(0)\hspace{1pt} :\ \sup_{\underline{\X} \in {\R}^{Nn}\by \cdots \by{\R}^{Nn^p}_s} \Phi(\underline{\X})
\big|\mF\big(x,u(x),\underline{\X}\big) \big|  \leq k \hspace{1pt} \right\}.
\]
Then, $\Om_{k}\sub \Om_{k+1}$ and $|\Om\set \Om_k| \larrow 0$ as $k\ri \infty$. We also define
\[
\Psi^k (x,\underline{\X})\hspace{1pt} :=\hspace{1pt} \Phi(\underline{\X})\big|\mF\big(x,u(x),\underline{\X}\big) \big| \chi_{\Om_k}(x), \quad k\in \N.
\]
Since $\de_{U^m}\weakstar \vartheta$ as $m\ri \infty$, we have 
\[
\int_{\Om} \Psi^k \big(x,U^m(x)\big)\hspace{1pt}dx \hspace{1pt} \larrow\hspace{1pt} \int_{\Om} \int_{\smash{\overline{\R}}^{Nn}\by \cdots \by\smash{\overline{\R}}^{Nn^p}_s}
\Psi^k \big(x,\underline{\X}\big)d[\vartheta(x)](\underline{\X})\hspace{1pt}dx.
\]
By assumption, we have $\mF\big(\cdot,u,[U^m]^R \big) \larrow 0$ a.e.\ on $\Om$ as $m\ri \infty$ and also the identity
\[
\Phi\left([U^m]^R\right) \mF\big(\cdot,u,[U^m]^R \big) \hspace{1pt} =\hspace{1pt} \Phi\left(U^m\right) \mF\big(\cdot,u, U^m\big)
\]
which is valid a.e.\ on $\Om$. From the above we infer that $\Psi^k (\cdot,U^m) \larrow 0$ a.e.\ on $\Om$. Moreover, by using the bound $|\Phi^k|\leq k$ and that $|\Om_k|<\infty$, the Dominated convergence theorem implies $\Psi^k (\cdot,U^m) \larrow 0$ in $L^1(\Om)$ as $m\ri \infty$. Hence, for a.e.\ $x\in \Om_k$ we have 
\[
\begin{split}
0\hspace{1pt}&=\hspace{1pt} \int_{\smash{\overline{\R}}^{Nn}\by \cdots \by\smash{\overline{\R}}^{Nn^p}_s}
\Psi^k\big(x,\underline{\X}\big)   \hspace{1pt} d[\vartheta(x)](\underline{\X}) \ \ \ \ \ \ \ \ \ 
\end{split}
\]
\[
\begin{split}
\ \ \ \ \ \  &=\hspace{1pt} \int_{\smash{\overline{\R}}^{Nn}\by \cdots \by\smash{\overline{\R}}^{Nn^p}_s}
\Phi(\underline{\X})\big|\mF\big(x,u(x),\underline{\X}\big) \big|  \hspace{1pt} d[\vartheta(x)](\underline{\X})
\\
&\geq\hspace{1pt} \int_{\mB_{R/2}(0)}
\big|\mF\big(x,u(x),\underline{\X}\big) \big|  \hspace{1pt} d[\vartheta(x)](\underline{\X}).
\end{split}
\]
The conclusion follows by letting $k\ri \infty$ and then $R\ri \infty$.

\ms
\noi (4)$\Rightarrow$(3): We argue as in the case ``(1)$\Rightarrow$(3)". Suppose that 
\[
\quad \quad \int_{ {{\R}}^{Nn} \by \cdots \by{{\R}}^{Nn^p}_s} 
\big| \mF\big(x,u(x),\underline{\X} \big) \big|\hspace{1pt} d[\vartheta (x) ](\underline{\X})\hspace{1pt} =\hspace{1pt} 0, \ \ \text{ a.e. }x\in \Om,
\]
whilst for some of these $x\in \Om$ we have $\supp\big(\vartheta_*(x)\big) \not\sub \big\{ \big|\mF\big(x,u(x),\cdot\big)\big|=0 \big\}$. Then, there exists $\underline{\X}_0$ with $\mF\big(x,u(x),\underline{\X}_0\big)=0$ with $[\vartheta(x)] \big(\mB_{R}(\underline{\X}_0) \big)>0$ for $R>0$. Thus, there exist $c_0,R_0>0$ such that $\left|\mF\big(x,u(x),\cdot\big)\right|\geq c_0>0$ on $\mB_{R_0}(\underline{\X}_0)$. Hence, $c_0\hspace{1pt} [\vartheta(x)] \big(\mB_{R_0}(\underline{\X}_0) \big)
 \leq 0$ and this contradiction establishes the desired inclusion.

\ms
\noi (3)$\Rightarrow$(6):  We fix $R>0$ and define $\Psi : \Om \by \smash{\overline{\R}}^{Nn} \by\cdots \by \smash{\overline{\R}}^{Nn^p}_s \hspace{1pt} \larrow\hspace{1pt} [0,\infty)$ by
\[
\Psi(x,\underline{\X})\hspace{1pt} :=\hspace{1pt}  \chi_{\overline{\mB_{R}(0)}} (\underline{\X})\hspace{1pt} \dist\left(\underline{\X}\ ,\hspace{1pt} {\mB_{R}(0)} \cap \Big\{\big|\mF\big(x,u(x),\cdot\big)\big|=0 \Big\} \right).
\]
Then, $\Psi$ is measurable in $x$ for all $\underline{\X}$ (this is a consequence of Aumann's theorem, see e.g.\ \cite{FL}), upper semicontinuous in $\underline{\X}$ for a.e.\ $x$ and also bounded. Hence, since $\de_{U^m}\weakstar \vartheta$ as $m\ri \infty$,  by Remark \ref{remark2}iii), we have 
\[
\begin{split}
\underset{m\ri \infty}{\lim \sup} & \int_{\Om}\Psi\big(x,U^m(x)\big) \hspace{1pt} dx\hspace{1pt} 
\leq\hspace{1pt} \int_{\Om} \int_{ \smash{\overline{\R}}^{Nn} \by\cdots \by \smash{\overline{\R}}^{Nn^p}_s } \Psi \big(x,\underline{\X}\big)\hspace{1pt} d[\vartheta(x)](\underline{\X})\hspace{1pt}dx\\
& =\hspace{1pt} \int_\Om \int_{\overline{\mB_{R}(0)}} 
\dist\left(\underline{\X}\ ,\ \Big\{\big|\mF\big(x,u(x),\cdot\big)\big|=0 \Big\}\cap \mB_{R}(0) \right) \hspace{1pt} d[\vartheta(x)](\underline{\X})\hspace{1pt}dx.
\end{split}
\]
By assumption and by Definition \ref{definition24a}, we have the inclusions $\supp\big(\vartheta(x)\big)  \cap \overline{\mB_{R}(0)} \sub \big\{\big|\mF\big(x,u(x),\cdot\big)\big|=0 \big\} \cap \overline{\mB_{R}(0)} \sub \{\Psi(x,\cdot)=0\}$, for a.e.\ $x\in\Om$. Hence, the last integral above vanishes and we obtain that $\Psi(\cdot,U^m) \larrow 0$ in $L^1(\Om)$ as $m\ri \infty$. Further, in view of Definition \ref{definition24}, we have the identity 
\[
\Psi\big(x,U^m(x)\big)\hspace{1pt} =\hspace{1pt}  \dist\left([U^m]^R(x)\ ,\  \Big\{\big|\mF\big(x,u(x),\cdot\big)\big|=0 \Big\} \cap \mB_{R}(0) \right),
\]
which holds for a.e.\ $x\in \Om$ and by using it we obtain that
\[
\int_\Om \dist\left([U^m]^R(x)\ ,\ \Big\{\big|\mF\big(x,u(x),\cdot\big)\big|=0 \Big\} \cap \mB_{R}(0) \right) \hspace{1pt}dx \hspace{1pt} \larrow \hspace{1pt} 0,
\]
as $m\ri \infty$. The conclusion follows by passing to a subsequence.   \qed
\ms

\subsection{Nonlinear nature of diffuse derivatives} (This subsection is not needed for  remainder of the paper.) In the context of classical PDE approaches (classical, strong, weak, distributional solutions), it is standard that the generalised derivative is a linear operator. However, this is generally \emph{false for diffuse derivatives}. \textit{Our approach is genuinely nonlinear and not a variant of classical developments}. Below we give a condition which guarantees that the sum of two $\mD$-solutions to a certain linear equation is a $\mD$-solution itself; this happens if at least one of the solutions is regular enough. Hence, the notions themselves are nonlinear even when we apply them to linear PDE. In order to proceed we need some notation. 

\begin{definition} \label{definition 32} Let $\mathbb{W}$ be a finite dimensional metric vector space isometrically and densely contained into a compactification $\mK$ of $\mathbb{W}$. Let also $T_a : \mW \ri \mW$ denote the translation operation given by $T_a b := b-a$. Given a probability $\vartheta \in \mP(\mK)$, we define $\vartheta \circ T_a \in \mP(\mK)$  by duality via the formula
\[
\ \ \ \langle \vartheta \circ T_a, \Phi \rangle\hspace{1pt} :=\hspace{1pt} \int_\mW \Phi(a+X)\hspace{1pt} d \vartheta (X)\hspace{1pt} + \int_{\mK\set\mW} \Phi(X)\hspace{1pt} d\vartheta (X),\ \ \ \ \Phi \in C(\mK).
\]
\end{definition}

Definition \ref{definition 32} requires translation of the part contained in the vector space while points ``at infinity" are left intact. 

\begin{proposition}[Diffuse derivatives \& $\mD$-solutions vs linearity] \label{proposition33} Let $u,v : \R^n \supseteq \Om \larrow \R^N$ be measurable maps. 

\noi a) If $v$ is differentiable in measure on $\Om$ with derivative $\mL \D v$, (Def.\ \ref{definition12}), then $\mD (u+v)= \mD u \circ T_{\mL \D v}$, a.e.\ on $\Om$. Here the diffuse Jets on both sides arise from the same infinitesimal sequence.

\noi b) Consider the measurable maps $\A^q : \R^n \supseteq \Om \larrow \R^{Nn^q}_s\!\otimes \R^M$ and $f,g :\R^n \supseteq \Om \larrow \R^M$ where $q=1,...,p$ and the linear systems $\A \!:: \! \D^{[p]}u = f$ and $\A \!:: \! \D^{[p]}v = g$.
Here $\A=(\A^1,...,\A^p)$. If $u,v$ are $\mD$-solutions, then $u+v$ is a $\mD$-solution $\A(x) \!:: \! \D^{[p]}(u+v) =f+g$, when $v$ is $p$-times differentiable in measure on $\Om$.
\end{proposition}
The notation ``::" above is a convenient abbreviation of the multiple contraction
\[
\sum_{\al_1, i_1}\A^1_{\mu ; \al_1, i_1} \D_{i_1}u_{\al_1}\hspace{1pt} +\hspace{1pt}  ... \hspace{1pt} +\sum_{\al_p, i^p_1 ...i^p_p } \A^p_{\mu ; \al_p ,i^1_1, ... ,i^p_p} \D^p_{i^p_1 ... i^p_p}u_{\al_p}  .
\]

The proof is based on the next general lemma.

\begin{lemma} \label{lemma34} Let $E\sub \R^n$ be measurable and $\mathbb{W}$ a finite dimensional metric vector space isometrically contained into a compactification $\mK$ of $\mathbb{W}$. If $U^m,V^m : E\sub \R^n \larrow \mW$ are measurable and such that $\de_{U^m}\weakstar \hspace{1pt}\vartheta$ in $\mY(E,\mK)$ and $V^m  \larrow\hspace{1pt} V$ a.e.\ on $E$, as $m\ri\infty$. Then, we have $\de_{U^m+V^m} \weakstar \hspace{1pt} \vartheta \circ T_V$ in $\mY(E,\mK)$ as $m\ri\infty$.
\end{lemma}

\BPL \ref{lemma34}. Fix $\phi \in L^1(E)$, $\Phi \in C(\mK)$ and $\e>0$. Since $\Phi$ is uniformly continuous, there is a \emph{bounded} increasing modulus of continuity $\om \in (C\cap L^\infty)[0,\infty)$ such that $|\Phi(X) -\Phi(Y)| \leq \om(|X-Y|)$ for $X,Y \in \mK$. Since $V^m\larrow V$ a.e.\ on $E$, we obtain $V^m \larrow V$ $\mu$-a.e.\ on $E$ where $\mu$ is the finite measure $\mu(A):=\|\phi\|_{L^1(A\cap E)}$, $A\sub \R^n$. It follows that $V^m \larrow V$ in $\mu$-measure as well. Hence, 
\[
\begin{split}
\left| \int_E \phi \Big[  \Phi(U^m +V^m) - \Phi(U^m+V)\Big] \right|\hspace{1pt} & \leq\hspace{1pt} \int_E |\phi|\hspace{1pt} \om \big(|V^m-V|\big)\\
 \leq \hspace{1pt} \|\om\|_{C(0,\infty)}  & \hspace{1pt}\mu\big(\{|V^m-V|>\e\}\big) \ +\ \om(\e)\hspace{1pt} \mu(E).
\end{split}
\]
By letting $m\ri \infty$ and then $\e\ri 0$, the density of the linear span of products  $\phi(x)\Phi(X)$ in $L^1\big(E,C(\mK)\big)$ and the definition of $\vartheta \circ T_V$ allow us to conclude.          \qed

\ms

\BPP \ref{proposition33}. If suffices to establish b) and only for $p=1$. By assumption, we have that $\A^1(x)\! :\! \mL \D v(x)=g(x)$ and also that for any $\Phi \in C_c^0(\R^{Nn})$, 
\[
\int_{\smash{\overline{\R}}^{Nn}} \Phi(X) \Big[\A^1(x): X-f(x)\Big]\hspace{1pt} d[\mD u(x)](X)\hspace{1pt}=\hspace{1pt}0,
\]
both being valid for a.e.\ on $x\in \Om$. Here $\mD u$ is any diffuse gradient. We fix any point $x$ as above and replace $\Phi$ by $\Phi\big(\cdot + \mL \D v(x)\big)$. Then, we obtain
\[
\int_{\smash{\overline{\R}}^{Nn}} \Phi \big(X + \mL \D v(x)\big) \Big[\A^1(x): \big(X + \mL \D v(x)\big)-f(x)-g(x)\Big]\hspace{1pt} d[\mD u(x)](X)\hspace{1pt}=\hspace{1pt}0.
\]
By the definition of $\mD u \circ T_{\mL \D v}$, we obtain
\[
\int_{\smash{\overline{\R}}^{Nn}} \Phi(Y) \Big[\A^1(x): Y-(f+g)(x)\Big]\hspace{1pt} d\big[\mD u(x) \circ T_{\mL \D v(x)}\big](Y)\hspace{1pt}=\hspace{1pt}0.
\]
By utilising part a), the conclusion ensues.   \qed

\begin{example}[Nonlinearity of diffuse derivatives] \label{example36} Let $K \sub \R$ be 
a compact nowhere dense set of positive measure (e.g. $K=[0,1]\set (\cup_1^\infty (r_j-3^{-j}, r_j+3^{-j}))$ where $(r_j)_1^\infty$ is an enumeration of $\mathbb{Q}\cap [0,1]$). Then, for $u:=\chi_K$ we have that $|\D^{1,h}u(x)|\ri \infty$ as $h \ri 0$ for $x\in K$ and $u'=0$ on $\R\set K$. Hence, by Lemma \ref{lemma2} along any $h_m\ri 0$ we have $\mD u(x)=\de_{\{ \infty \}}$ for a.e.\ $x\in K$. However, for $v:=-u$, we have $\mD (u+v)=\de_{\{0\}}$ a.e.\ on $\R$, while $\mD u = \mD v=\de_{\{ \infty \}}$ a.e.\ on $K$.
\end{example}

\noi \textbf{Comparison with distributional solutions}. Let us conclude this section with an \emph{informal} discussion of the relation between distributional and $\mD$-solutions. Let us first compare distributional to diffuse derivatives. First recall that the distributional gradient $\D u$ of $u\in L^1_{\text{loc}}(\R^n)$ can be weakly* approximated by difference quotients: for any $\phi\in C^\infty_c(\R^n)$,
\[ 
\langle \phi, \D u \rangle\hspace{1pt} =\hspace{1pt} \lim_{m\ri\infty} \int_{\R^n}\phi\hspace{1pt} \D^{1,h_{m}}u\hspace{1pt} =\hspace{1pt} \lim_{m\ri\infty} \int_{\R^n}\phi\left( \int_{\R^{n}} X\hspace{1pt} d \big[\de_{\D^{1,h_{m}}u} \big](X)\right).
\]
If ``bar$_*$" denotes the barycentre of the restriction of a measure on $\smash{\overline{\R}}^n$ off $\{\infty\}$, the above can be rewritten as
\beq \label{2.14}
\text{bar}_*\left(\de_{\D^{1,h_{m}}u} \right) \hspace{1pt}\weakstar \hspace{1pt} \D u, \quad \text{ as }m \ri\infty,
\eeq
in the distributions $\mathscr{D}'\big(\R^n,\R^{n}\big)$. Along perhaps a further subsequence, we have
\beq \label{2.15}
{\de_{\D^{1,h_{m }}u}} \hspace{1pt}\weakstar \hspace{1pt} \mD u,\ \text{ in }\mY\big(\R^n,\smash{\overline{\R}}^n\big), \quad \text{ as }m \ri\infty.
\eeq
By juxtaposing \eqref{2.14} with \eqref{2.15}, our interpretation is that \textit{the barycentre of the diffuse derivative (off $\{\infty\}$) is unique and equals the distributional derivative}: $ \text{bar}_*( \mD u ) = \D u$. Regarding the notions to solution, apparently $\mD$-solutions are a more general theory than distributional solutions in the sense that they apply to more general PDEs and under weaker requirements. However, the two theories are not immediately comparable on their common domain of $L^1_{\text{loc}}$ solutions of linear systems with smooth coefficients. On the one hand, Proposition \ref{proposition33} and Example \ref{example36} point out a property which is not generally true for diffuse derivatives but is always true for distributional derivatives. However, $\mD$-solutions completely avoid the impossibility to multiply distributions. For example, if $\A \in L^\infty(\R^n,\R^n)$,
\[
 {\A \cdot \D^{1,h_{m }}u \hspace{1pt} -\!\!\!- \! \Big\langle}
\!\!\!{\begin{split}
 &\ \ \ { /\!\!\!\!\!\!\!\!\!\weakstar \hspace{1pt} \A \cdot \D u, \ \text{ in } \mathscr{D}'\big(\R^n,\R^{n}\big), \text{ \ \ [not well defined!]}} \phantom{\big|_{j}}\\
& {\weakstar \A \cdot \mD u, \ \hspace{1pt} \text{ in } \mY\big(\R^n,\smash{\overline{\R}}^n\big).  \text{ \ \ \hspace{1pt} [well defined!]}}
\end{split}}
\]
Hence, although the product $\A(x)\cdot \D u(x) = \A\cdot  \text{bar}_*\left(\mD u(x) \right)$ is ill-defined, diffuse derivatives make sense because they can be multiplied with measurable functions.

\section{$\mD$-solutions of the $\infty$-Laplacian and tangent systems} \label{section3}

In this section we establish our first main result concerning $\mD$-solutions. We treat the Dirichlet problem for the $\infty$-Laplace system \eqref{1.9} which is the fundamental equation of vectorial Calculus of Variations in the space $L^\infty$ and arises from the functional \eqref{1.8}.

\begin{theorem}[Existence of $\infty$-Harmonic maps] \label{theorem20}  Let $\Om \sub \R^n$ be an open set with $|\Om|<\infty$, $n\geq 1$. For any $g\in W^{1,\infty}(\Om,\R^n)$, the Dirichlet problem
\beq \label{3.1}
\left\{
\begin{array}{rl}
\ \ \De_\infty u \hspace{1pt} =\hspace{1pt} 0, & \text{ in }\Om,\\
u\hspace{1pt} =\hspace{1pt} g, & \text{ on }\p\Om,
 \end{array}
 \right.
\eeq
has a $\mD$-solution $u \in W_g^{1,\infty}(\Om,\R^n)$ with respect to the \emph{standard frames} (Definition \ref{definition11}). In particular, for any $\mD^2u \in \mY\big(\Om, \smash{\overline{\R}^{nn^2}_s} \big)$ and $\Phi\in C_c\big( \smash{\R^{nn^2}_s}\big)$, we have
\[
\ \  \ \ \ \int_{\smash{\overline{\R}}^{nn^2}_s} \Phi(\X) \hspace{1pt} \Big(\D u \ot \D u + |\D u|^2[\![ \D u ]\!]^\bot \! \ot \mathrm{I} \Big):\X\hspace{1pt} d[\mD ^2 u](\X)\hspace{1pt} =\hspace{1pt} 0, \ \ \text{ a.e.\ on $\Om$.}
\]
\end{theorem}

\begin{remark} \label{rem30} Unfortunately, as we proved in \cite{K2}, in general it is impossible to obtain uniqueness of solutions to the equations of vectorial $L^\infty$ problems \textbf{even within the class of smooth solutions}. Clearly, uniqueness in the vectorial case is not an issue of defining a ``proper" notion of generalised solution, since even classical solutions in general are non-unique. Instead, \emph{extra conditions} must be determined to select a ``good" solution. On the other hand, uniqueness is standard in the scalar case (a celebrated theorem of Jensen, see e.g.\ \cite{C}, \cite{K8}). Such phenomena are \emph{not exclusive to the $\infty$-Laplacian}: for instance, the Dirichlet problem for the minimal surface system may have either non-existence or non-uniqueness in codimension greater than one (see \cite{OL}), while for the minimal surface equation it is well posed. 
\end{remark}

In addition, the next corollary will also be established in the course of its proof.

\begin{corollary}[Multiplicity \& geometric properties of $\mD$-solutions] \label{corollary21} In the setting of Theorem \ref{theorem20}, if $n\geq 2$ then \eqref{3.1} actually has an infinite set of solutions. Moreover, for any $M> \big\| (\D g^\top \! \D g)^{1/2}\big\|_{L^\infty(\Om)}$ there is a $\mD$-solution $u=u(M)$ satisfying
\beq \label{3.4}
|\D u|^2\hspace{1pt} =\hspace{1pt} nM^{2}, \ \ \ \big|\!\det(\D u)\big|\hspace{1pt} =\hspace{1pt} M^{n}, \ \ \text{ a.e.\ on }\Om.
\eeq
\end{corollary}

Hence, the $\mD$-solutions we construct are ``critical points" with pre-assigned energy level, having also the geometric property of being full-rank solutions of the vectorial Eikonal equation. 

\subsection{The idea of the proof} \label{subsection3.1} Suppose that $u \in C^2(\Om,\R^n)$ solves \eqref{1.9} and recall that $[\![ \D u ]\!]^\bot= \text{Proj}_{(R(\D u))^\bot}$. By contracting derivatives, we rewrite the system as
\beq \label{3.5}
\D u\hspace{1pt}\D\Big(\frac{1}{2}|\D u|^2 \Big)\hspace{1pt} +\hspace{1pt} |\D u|^2[\![ \D u ]\!]^\bot \De u\hspace{1pt} =\hspace{1pt} 0.
\eeq
It follows that \emph{smooth} solutions of the $1$st order differential inclusion $\D u(x) \in \mathcal{K}_c$, $x\in \Om$, where $c>0$ is a parameter and 
\[
\mathcal{K}_c\hspace{1pt}:=\hspace{1pt} \Big\{ X\in \R^{nn} \hspace{1pt} : \hspace{1pt} |X|\hspace{1pt} =\hspace{1pt} c, \hspace{1pt}  |\!\det(X)| >\hspace{1pt} 0 \Big\} , 
\]
actually are $\infty$-Harmonic mappings: indeed, if $\D u(\Om)\sub \mathcal{K}_c$, then $|\D u|^2\equiv c^2$ and $\det(\D u)\neq 0$ on $\Om$. In view of \eqref{3.5} we have that the system is satisfied because $|\D u| \equiv \text{const}$ and $u$ is a submersion which gives $[\![ \D u ]\!]^\bot\equiv 0$. Hence, if we prove existence of a solution to the inclusion with the desired boundary data, this yields a solution to \eqref{3.1}. However, the preceding arguments make sense only for classical or strong solutions. The starting point of the proof of  Theorem \ref{theorem20} is to use the Dacorogna-Marcellini Baire Category method \cite{DM} in order to construct Lipschitz solutions of the inclusion with the given boundary data. Then, by using the machinery of $\mD$-solutions we make the previous ideas rigorous for Lipschitz maps, which is the natural regularity class. Note that our methodology is \emph{not variational} and does not directly involve the functional \eqref{1.8}.

\subsection{Proof of the main result} The case $n=1$ is trivial (see \cite{K1}), so we may assume $n\geq 2$. A central ingredient in the proof of Theorem \ref{theorem20} is a result of independent interest, Theorem \ref{lemma31} below, which provides a method of constructing nonsmooth $\mD$-solutions to nonlinear systems by ``differentiating an equation".

\begin{theorem}[Differentiating equations in the $\mD$-sense]  \label{lemma31} 
Let $\mF$ be a $C^1$ map as in \eqref{1.1} and consider the $p$-th order system
\[
\mF\left(x,u(x),\D^{[p]}u(x)\right)\hspace{1pt}=\hspace{1pt} 0, \ \ \ x\in \Om.
\]
If  $u \in W^{p,\infty}_{\text{loc}}(\Om,\R^N)$ is a strong a.e.\ solution to the system, then $u$ is a $\mD$-solution to the ``tangent system" on $\Om$ with respect to the usual frames (Definition \ref{definition13}):
\[
\mF_x\big(\cdot,u,\D^{[p]}u\big)\hspace{1pt} +\hspace{1pt} \mF_\eta \big(\cdot,u,\D^{[p]}u\big)\D u \hspace{1pt}+\hspace{1pt} \mF_{\underline{\X}} \big(\cdot,u,\D^{[p]}u\big)\!::\!\D^{[p+1]}u\ =\ 0.
\]
\end{theorem}

For the notation ``::" see Proposition \ref{proposition33}. Theorem \ref{lemma31} is actually true for solutions which are merely $W^{p,1}_{\text{loc}}(\Om,\R^N)$ or just $p$-times differentiable in measure (Definition \ref{definition12}), but then we have to assume certain growth bounds on the derivatives of $\mF$. We invite the reader to note the \emph{simplicity} with which we pass to limits in the proof below within the framework of $\mD$-solutions.

\BPT \ref{lemma31}. It suffices to prove only the case of $p=1$ and with no explicit $u$ dependence, the general case following analogously. Hence we suppose that $u\in W^{1,\infty}_{\text{loc}}(\Om,\R^N)$ solves $\mF\big(x,\D u(x)\big)=0$ for a.e.\ $x\in \Om$ and we aim at showing
\[
\mF_{x} \big(x,\D u(x)\big)\hspace{1pt} +\hspace{1pt} \mF_X  \big(x,  \D u(x)   \big)  : \D^2 u(x)  \hspace{1pt} =\hspace{1pt} 0, \ \ \ x \in \Om,
\]
in the $\mD$-sense (Definition \ref{definition11}). For a.e.\ point $x\in \Om$ such that $\mF\big(x,\D u(x)\big)=0$ and $h\neq 0$ small enough, Taylor's theorem implies for each $i$ the identity
\beq \label{id}
\begin{split}
\mF_{x_i} &\big(x,\D u(x)\big)\hspace{1pt} +\hspace{1pt} \mF_X  \big(x,  \D u(x)   \big)  : \D^{1,h}_{e^i}\D u(x)  \\
 = & \ -  \D^{1,h}_{e^i}\D u(x) : \int_0^1 \Big\{\mF_X\Big(x+\la he^i, \D u(x)+\la \big[ \D u(x+he^i)-\D u(x)\big]\Big) \\
 &\hspace{90pt}  - \mF_X\big(x,\D u(x)\big) \Big\} \hspace{1pt} d\la \\
&\ - \int_0^1 \Big\{\mF_{x_i}\Big(x+\la he^i, \D u(x) + \la \big[ \D u(x+he^i)-\D u(x)\big]\Big)\\
 &\hspace{40pt}   - \mF_{x_i}\big(x,\D u(x)\big) \Big\} \hspace{1pt} d\la.
\end{split}
\eeq
We fix an infinitesimal sequence $(h_m)_{m=1}^\infty\sub \R\set\{0\}$ and observe that by the weak* compactness of $\mY\big(\Om,\smash{\overline{\R}}^{Nn^2}_s \big)$, along a subsequence $h_{m_k}\ri 0$ we have
\[
\de_{\D^{1,h_{m_k}}\D u}\weakstar \mD ^2 u \ \hspace{1pt} \text{ in }\ \mY\big(\Om,\smash{\overline{\R}}^{Nn^2}_s \big),\ \hspace{1pt} \text{ as } k\ri \infty.
\]
We now invoke that $|\D u| \in L^\infty(\Om)$ to infer that since $\D u(\cdot+he^i)\larrow \D u$ in $L^1(\Om,\R^{Nn})$ as $h\ri 0$, there is a further subsequence denoted again by $(h_{m_{k}})_{k=1}^\infty$ such that for a.e.\ $x\in \Om$ we have $\D u(x+h_{m_{k}}e^i)\larrow \D u(x)$ as $k\ri \infty$. Next, we set
\[
G^\infty_i(x,\X)\ := \  \mF_{x_i} \big(x,\D u(x)\big) \hspace{1pt} +\ \sum_{\be,j} \mF_{X_{\be j}}  \big(x,  \D u(x)   \big) \hspace{1pt} \X_{\be ji}
\]
and for $m\in \N$
\[
\begin{split}
G^m_i(x,\X)\hspace{1pt} &:=  \ \mF_{x_i} \big(x,\D u(x)\big) \hspace{1pt} +\ \sum_{\be,j} \mF_{X_{\be j}}  \big(x,  \D u(x)   \big)  \hspace{1pt} \X_{\be ji} \\
&\hspace{1pt}  +  \sum_{\be,j} \X_{\be ji} \int_0^1 \Big\{\mF_{X_{\be j}} \Big(x+\la h_me^i, \D u(x)+\la \big[ \D u(x+h_me^i)-\D u(x)\big]\Big) \\
 &\hspace{70pt}  - \mF_{X_{\be j}} \big(x,\D u(x)\big) \Big\} \hspace{1pt} d\la \\
& + \int_0^1 \Big\{\mF_{x_i}\Big(x+\la h_me^i, \D u(x) + \la \big[ \D u(x+h_me^i)-\D u(x)\big]\Big) \\
 &\hspace{30pt}   - \mF_{x_i}\big(x,\D u(x)\big) \Big\} \hspace{1pt} d\la.
\end{split}
\]
By the $C^1$ regularity of $\mF$ and that $\D u(\cdot+h_{m_k}e^i)\larrow \D u$ a.e.\ on $\Om$ as $k\ri \infty$ (together with the Dominated convergence theorem and that $\D u \in L^\infty_{\text{loc}}(\Om,\R^{Nn})$), for a.e.\ $x\in \Om$ we obtain that $G^{m_k} (x,\cdot)\larrow G^\infty(x,\cdot)$ in $C\big( \R^{Nn^2}_s,\R^M \big)$, as $k\ri \infty$. Moreover, in view of the definition of $G^m$, the identity \eqref{id} gives
\[
G^m\Big(x,\D^{1,h_m}\D u(x) \Big)\hspace{1pt}=\hspace{1pt} 0\hspace{1pt} \ \text{ a.e.\ on }\Om,\ \hspace{1pt} m\in\N.
\]
Hence, for any $\Phi\in C_c\big( \R^{Nn^2}_s \big)$ we have
\[
\int_{\smash{\overline{\R}}^{Nn^2}_s} \Phi(\X)\hspace{1pt} G^{m_k}(x,\X)\hspace{1pt} d\big[\de_{\D^{1,h_{m_k}}\D u(x)} \big] (\X)\hspace{1pt} =\hspace{1pt} 0 \ \ \ \text{a.e. }x\in\Om,
\]
for $k\in \N$. The convergence Lemma \ref{lemma16A} now implies 
\[
\int_{\smash{\overline{\R}}^{Nn^2}_s} \Phi(\X)\hspace{1pt} G^\infty(x,\X)\hspace{1pt} d\big[ \mD^2 u(x) \big] (\X)\hspace{1pt} =\hspace{1pt} 0, \ \ \ \text{a.e. }x\in\Om,
\]
for any $\Phi\in C_c\big( \R^{Nn^2}_s \big)$ and any diffuse hessian $\mD ^2 u \in \mY\big(\Om,\smash{\overline{\R}}^{Nn^2}_s\big)$. Hence, $u$ is a $\mD$-solution of $G^\infty\big(\cdot,\D^2 u \big)=0$ and by the definition of $G^\infty$, the result ensues.        \qed

\ms

\BPT \ref{theorem20} (and Corollary \ref{corollary21}). Assume we are given $\Om \sub \R^n$ with finite measure and $g\in W^{1,\infty}(\Om,\R^n)$. We begin with the next:

\bc \label{claim22} If $M> \| (\D g^\top \D g)^{1/2}\|_{L^\infty(\Om)}$, there exists $u\in W^{1,\infty}_g(\Om,\R^n)$ such that $|\D u|^2= nM^{2}$ and also $|\det(\D u)|= M^{n}$, both holding a.e.\ on $\Om$.
\ec

\BPC  \ref{claim22}. Given a map $u: \R^n \supseteq \Om \larrow \R^n$ in $W^{1,\infty}_g(\Om,\R^n)$, let $\la_i(\D u)$ denote the $i$th singular value, that is the $i$th eigenvalue of $(\D u^\top \D u)^{1/2}$:
\[
\si\big( (\D u^\top \D u)^{1/2} \big)\hspace{1pt} =\hspace{1pt} \big\{\la_1(\D u),\ldots,\la_n(\D u) \big\}, \ \ \la_i\leq\la_{i+1}.
\] 
Fix an $M>0$ as in statement and consider the Dirichlet problem:
\beq \label{DP}
\left\{
\begin{array}{rl}
\la_i(\D v)\hspace{1pt} =\hspace{1pt} 1, &\ \ \text{ a.e. in }\Om,\ \ i=1,...,n,\\
v\hspace{1pt} =\hspace{1pt} g/M, & \ \ \text{ on }\p \Om.
\end{array}
\right.
\eeq
Then, we have the estimate
\beq \label{CC}
\ \ \ \big\| \la_n(\D g)\big\|_{L^\infty(\Om)} =\hspace{1pt} \Big\| \max_{|e|=1}\hspace{1pt} (\D g^\top \D g)^{1/2}:e\ot e\Big\|_{L^\infty(\Om)} \hspace{1pt} \leq \hspace{1pt} \big\| (\D g^\top \D g)^{1/2}\big\|_{L^\infty(\Om)}.
\eeq
In view of the results of \cite{DM}, the estimate \eqref{CC} implies that the required compatibility condition is satisfied in regard to the problem \eqref{DP}. Hence there is a strong solution $v$ to \eqref{DP} such that $v-(g/M) \in W_0^{1,\infty}(\Om,\R^n)$ for the given $M$ and the boundary data $g$. Finally, since $\la_i(\D v)=1$ a.e.\ on $\Om$, by setting $u:=Mv$ we have
\[
\begin{split}
|\D u|^2&=\hspace{1pt}  M^2|\D v|^2\hspace{1pt} =\hspace{1pt}  M^2 \sum_{i=1...n} \la_i(\D v)^2\hspace{1pt} =\hspace{1pt} nM^2,\ \ \text{ a.e.\ on }\Om,\\
|\det(\D u)|\hspace{1pt} &=\hspace{1pt}  M^n |\det(\D v)| \hspace{1pt} =\hspace{1pt}   M^n \prod_{i=1...n} \la_i(\D v)\hspace{1pt} =\hspace{1pt} M^n, \ \ \text{ a.e.\ on }\Om,
\end{split}
\]
and in addition, $u-g\in W_0^{1,\infty}(\Om,\R^n)$. The proof of the claim is complete.      \qed

\smallskip

Now we may complete the proof. For the given boundary condition $g$, we fix an $M>0$ as in the claim and consider one of its solutions $u\in W^{1,\infty}_g(\Om,\R^N)$ which satisfies $|\D u|^2-nM^2=0$, a.e.\ on $\Om$. We set $\mF(X) :=  |X|^2-nM^2$ for $X\in \R^{Nn}$ and apply Theorem \ref{lemma31} to infer that $u\in W^{1,\infty}_g(\Om,\R^N)$ is a $\mD$-solution to the tangent system $\mF_X(\D u):\D^2 u=0$; that is, for all $i$ we have
\[
\sum_{\be, j} \D_j u_\be (x)\hspace{1pt} \D^2_{ij}u_\be(x)\hspace{1pt} =\hspace{1pt} 0, \ \ \ x\in \Om, \text{ in the $\mD$-sense}.
\]
This means that when $\de_{\D^{1,h_{m}}\D u}\weakstar \mD ^2 u$ in $\mY\big(\Om,\smash{\overline{\R}}^{nn^2}_s \big)$ as $m\ri \infty$, we have
\[
\int_{\smash{\overline{\R}}^{nn^2}_s} \sum_{\be, j} \D_j u_\be (x)\hspace{1pt} \Phi(\X) \hspace{1pt} \X_{\be i j} \hspace{1pt} d\big[\mD ^2 u(x)\big](\X)\hspace{1pt} =\hspace{1pt} 0, \ \ \ \text{ a.e.\ }x\in\Om,
\]
for any fixed $\Phi\in C_c\big(\R^{nn^2}_s\big)$. We multiply the above equation by $\D_i u_\al(x)$ and sum with respect to $i$ to obtain
\[
\int_{\smash{\overline{\R}}^{nn^2}_s} \sum_{\be, j,i} \Phi(\X) \hspace{1pt} \D_i u_\al(x)\hspace{1pt} \D_j u_\be (x)\hspace{1pt} \X_{\be i j} \hspace{1pt} d\big[\mD ^2 u(x)\big](\X)\hspace{1pt} =\hspace{1pt} 0, \ \ \ \text{ a.e.\ }x\in\Om.
\]
Finally, by Claim \ref{claim22} we have $\det(\D u)\neq 0$ a.e.\ on $\Om$ and as a result $\D u(x)$ has rank equal to $n$ in $\R^{nn}$. Hence, the projection $[\![ \D u(x)]\!]^\bot=\text{Proj}_{R(\D u(x))^\bot}$ vanishes for a.e.\ $x\in \Om$ and consequently we obtain
\[
\int_{\smash{\overline{\R}}^{n n^2}_s} \sum_{\be, i} \Phi(\X) \hspace{1pt} \big|\D u(x)\big|^2 [\![\D u(x)]\!]^\bot_{\al \be} \hspace{1pt} \X_{\be i i} \hspace{1pt} d\big[\mD ^2 u(x)\big](\X)\hspace{1pt} =\hspace{1pt} 0, \ \ \ \text{ a.e.\ }x\in\Om,
\]
for any $\Phi\in C_c\big(\R^{nn^2}_s\big)$ and any diffuse hessian $\mD ^2 u \in \mY\big(\Om,\smash{\overline{\R}}^{n n^2}_s \big)$. The last two equations imply that $u$ is $\infty$-Harmonic in the $\mD$-sense and the theorem follows.          \qed

\ms

\section{$\mD$-solutions of fully nonlinear degenerate elliptic systems} \label{section4}

Fix $n,N\geq 1$, let $\Om \sub \R^n$ be an open set and $\mF : \Om \by \R^{Nn^2}_s \larrow \R^N$ a Carath\'eodory map. In this section we establish our second main result, namely the \emph{existence of a unique $\mD$-solution} $u:\R^n \supseteq \Om \larrow \R^N$ to the Dirichlet problem
\beq \label{4.1}
\left\{
\begin{split}
\mF(\cdot,\D^2 u)\hspace{1pt} &=\hspace{1pt} f, \ \ \text{ in }\Om,\\
u\hspace{1pt}&=\hspace{1pt} 0,\ \ \text{ on }\p\Om, 
\end{split}
\right.
\eeq
together with a \emph{partial regularity assertion of new type involving differentiability along rank-one directions} (instead of the usual partial regularity on a subset of the domain). We will assume that $f\in L^2(\Om,\R^N)$ and $\mF$ satisfies a degenerate ellipticity condition which in general does not guarantee that solutions are {even \emph{once} weakly differentiable}. This extends previous results of the author in the class of strong solution for \eqref{4.1} (\cite{K9,K11}) under a stronger ellipticity condition.

\subsection{The idea of the proof} The solvability of \eqref{4.1} in the class of $\mD$-solutions is based on  the study of the linearised system with constant coefficients
\beq \label{4.2}
\left\{
\begin{split}
\A : \D^2 u\hspace{1pt} &=\hspace{1pt} f, \ \ \text{ in }\Om,\\
u\hspace{1pt}&=\hspace{1pt} 0, \ \ \text{ on }\p\Om, 
\end{split}
\right.
\eeq
when $\A $ is a (perhaps degenerate) symmetric $4$th order tensor and on a perturbation device provided by our ellipticity assumption for $\mF$. The latter allows to solve \eqref{4.1} by solving \eqref{4.2} and using a fixed point argument in the guises of a classical theorem of Campanato (\cite{C3}).  In order to solve \eqref{4.2} in the $\mD$-sense (and not just weakly) we impose a structural condition on $\A$ which allows to construct $\mD$-solutions as maps having weakly twice differentiable projections along certain rank-one lines of $\R^{Nn}$. These are the ``directions of ellipticity" of \eqref{4.2}. We formalise this idea by introducing a ``fibre" extension of the classical Sobolev spaces which consist of maps possessing only certain partial regularity along rank-one lines. Our fibre space counterparts are adapted to the degenerate nature of the problem and support feeble yet sufficient versions of weak compactness, trace operators and Poincar\'e inequalities established in \cite{K10}. The proof is completed by characterising the fixed point as the unique $\mD$-solution of the problem \eqref{4.1} in the fibre space.

\subsection{Fibre spaces, degenerate ellipticity and the main result} \label{subsection4.2} Before stating our existence result we need some preparation. We will use the notation $\A \in \R^{Nn \by Nn}_s$ to symbolise symmetric linear maps $\A :\R^{Nn}\larrow \R^{Nn}$, i.e.\ 4th order tensors satisfying $\A_{\al i \be j }=\A_{\be j \al i}$ for all $\al,\be =1,...,N$, $i,j =1,...,n$. The notation $N\big(\A : \R^{Nn}\ri \R^{Nn} \big)$ and $N\big(\A : \R^{Nn^2}_s\ri \R^{N} \big)$ will be used to symbolise the nullspaces of $\A$ when it acts as a linear map with domain and range like those indicated in the brackets, i.e. $Q \lmapsto \sum_{\al,\be,i,j}\big(\A_{\al i \be j }Q_{\be j}\big)\hspace{1pt} e^\al \ot e^i$ and $\X \lmapsto \sum_{\al,\be,i,j}\big( \A_{\al i \be j }\X_{\be i j}\big)\hspace{1pt} e^\al$. We will also use similar notation for the respective ranges with ``$R$" instead of ``$N$". If $\A$ is rank-one non-negative, i.e.\ if the respective quadratic form is rank-one convex $\A : \eta \ot a \ot \eta \ot a\hspace{1pt} = \sum_{\al,\be,i,j}\A_{\al i \be j }\eta_\al \hspace{1pt} a_i\hspace{1pt} \eta_\be\hspace{1pt} a_j\geq 0$ ($\eta \in \R^N,\ a\in \R^n$), we define 
\beq \label{4.3}
\begin{split}
 \Pi \hspace{1pt} &:=\ R\Big(\A : \R^{Nn}\ri \R^{Nn} \Big)  \hspace{78pt}    \sub\  \R^{Nn},\\
 \Si \hspace{1pt} &:=\ \spn [\Big\{ \hspace{1pt}\eta \ \Big|\  \eta \ot a \in \Pi\hspace{1pt} \Big\}] \hspace{60pt} \ \ \ \sub\ \R^N,\\
 \Xi \hspace{1pt} &:=\ \spn [\Big\{ \hspace{1pt}\eta \ot (a\vee b)  \ \Big|\ \eta \ot a,\ \eta \ot b \in \Pi\hspace{1pt} \Big\}] \sub\  \R^{Nn^2}_s, \\
 \nu \hspace{1pt}&:=\ \min_{|\eta|=|a|=1, \hspace{1pt}\eta \ot a\in \Pi} \Big\{ \A : \eta \ot a \ot \eta \ot a \Big\}\ \hspace{1pt}\hspace{1pt}\hspace{1pt} \ \ \ >\hspace{1pt} 0.
\end{split}
\eeq
We will call $\nu$ the \emph{ellipticity constant of $\A$}, bearing in mind that strictly speaking $\A$ may not be elliptic and the respective infimum over $\R^{Nn}$ may vanish. We also recall that we will use the same letters $\Pi,\Xi,\Si$ to symbolise the subspaces as well as the orthogonal projections on them.

\smallskip

\noi \textbf{The fibre Sobolev spaces.}  Given $\A\in\R^{Nn \by Nn}_s$ rank-one non-negative, let $\Si, \Pi,\Xi$ be as in \eqref{4.3} and suppose \emph{$\Pi$ is spanned by rank-one directions}. A sufficient condition for this to happen is when $\A$ is decomposable (Definition \ref{SH} that follows). For simplicity, we treat only the $L^2$ second order case needed in this paper. We begin by identifying $W^{2,2}(\Om,\R^N)$ with its isometric image $\tilde{W}^{2,2}(\Om,\R^N)$ into a product of $L^2$ spaces $\tilde{W}^{2,2}(\Om,\R^N) \underset{^\ri}{\subset} L^2\big(\Om,\R^N \by \R^{Nn} \by \R^{Nn^2}_s\big)$ via the map $u\lmapsto(u,\D u,\D^2 u)$. We define the \textbf{fibre Sobolev space} $\mathscr{W}^{2,2}(\Om,\Si)$ as the Hilbert space
\beq \label{4.4A}
\mathscr{W}^{2,2}(\Om,\Si)\hspace{1pt} :=\hspace{1pt} \overline{\hspace{1pt} \text{Proj}_{ L^2\big(\Om,\Si \by \Pi \by \Xi\big) }\hspace{1pt}  \tilde{W}^{2,2}(\Om,\R^N) \hspace{1pt} }^{\|\cdot\|_{ L^2(\Om) } }
\eeq
which we equip with the norm (written for $W^{2,2}$ maps)
\[
\|u\|_{\mathscr{W}^{2,2}(\Om,\Si)} \hspace{1pt} :=\ \big\|\Si \hspace{1pt} u \big\|_{L^2(\Om)}+\hspace{1pt} \big\|\Pi \hspace{1pt} D u \big\|_{L^2(\Om)}+\hspace{1pt} \big\|\Xi\hspace{1pt} \D^2 u \big\|_{L^2(\Om)}.
\]
By the Mazur theorem,  $\mathscr{W}^{2,2}(\Om,\Si)$ can be characterised as
\[
\mathscr{W}^{2,2}(\Om,\Si)\hspace{1pt} = \left\{ 
\begin{array}{l}
 \left(u,  G  (u),\mathrm{G}^2(u)\right) \in\hspace{1pt} L^2\big(\Om,\Si \by \Pi \by \Xi\big) \ \big|\  \exists\ (u^m)_1^\infty   \sub \ms \\
W^{2,2}(\Om,\R^N) :\ \hspace{1pt}  \text{ we have weakly in }L^2 \text{ as }m\ri\infty \ms\\
\text{that }\hspace{1pt} \left(\Si\hspace{1pt}u^m , \Pi\hspace{1pt}\D u^m, \Xi \hspace{1pt}\D^2 u^m \right) \weak \big(u, \mathrm{G}(u) , \mathrm{G}^2(u)\big)
\end{array}
\right\}.
\]
We will call $\mathrm{G}(u) \in L^2(\Om,\Pi)$ the \textbf{fibre gradient} of $u$ and $\mathrm{G}^2(u)\in L^2(\Om,\Xi)$ the \textbf{fibre hessian} of $u$. 

We now show that $\big(\mathrm{G}(u),\mathrm{G}^2(u)\big)$ \emph{depend only on $u \in L^2(\Om,\Si)$ and not on the approximating sequence}. Indeed, let $(u^m)_1^\infty$ and $(v^m)_1^\infty$ be sequences in ${W}^{2,2}(\Om,\R^N)$ such that
\[
\begin{split}
\left(\Si\hspace{1pt}u^m , \Pi\hspace{1pt}\D u^m, \Xi \hspace{1pt}\D^2 u^m \right) \weak \big(u, \mathrm{G}(u), \mathrm{G}^2(u)\big),
\\
\left(\Si\hspace{1pt}v^m , \Pi\hspace{1pt}\D v^m, \Xi \hspace{1pt}\D^2 v^m \right) \weak \big(v, \mathrm{G}(v), \mathrm{G}^2(v)\big),
\end{split}
\]
weakly in $L^2\big(\Om,\Si \by\Pi\by\Xi\big)$ as $m\ri\infty$. We immediately have that $\Si u =\Si v$ a.e.\ on $\Om$ and hence $u, v$ represent the same element of $L^2(\Om,\Si)$ because their projections on the subspace $\Si\sub \R^N$ coincide, whilst by definition $\Si^\bot u \equiv \Si^\bot v \equiv 0$. Similarly, since $\Si,\Pi,\Xi$ are spanned by directions of the form $\eta$, $\eta \ot a$ and $\eta \ot (a\vee b) $ respectively, for any $\phi \in C^\infty_c(\Om)$, $\eta \ot a \in \Pi$ and $\eta \in \Si$ we have
\[
\begin{split}
& \int_\Om \phi \hspace{1pt}\big(G(u)-G(v)\big):\eta \ot a \hspace{1pt} =\hspace{1pt} \lim_{m\ri \infty} \int_\Om \phi \hspace{1pt} \big(\Pi\hspace{1pt} \D u^m-\Pi\hspace{1pt} \D v^m\big):\eta \ot a
\\
&=\hspace{1pt} \lim_{m\ri \infty} \int_\Om \phi \hspace{1pt}\D_a\left[\eta \cdot \big(u^m -v^m\big)\right]
\hspace{1pt}=\hspace{1pt} -\lim_{m\ri \infty} \int_\Om \D_a\phi \hspace{1pt}\left[\eta \cdot \Si\big(u^m -v^m\big)\right] \hspace{1pt} =\hspace{1pt} 0
\end{split}
\]
and hence $G(u), G(v)$ coincide as elements of $L^2(\Om,\Pi)$ because their projections on the subspace $\Pi \sub \R^{Nn}$ coincide. The remaining case is analogous.

Further, by using the standard properties of equivalence between strong and weak $L^2$ directional derivatives, we have that $\mathrm{G}(u),\mathrm{G}^2(u)$ can be characterised as ``fibre" derivatives of $u$: for any directions $\eta \in \Si$, $\eta \ot a \in \Pi$ and $\eta \ot (a\vee b) \in \Xi$, we have $\mathrm{G}(u) : (\eta \ot a) = D_a(\eta \cdot u)$ and also
\[
\begin{split}
\mathrm{G}^2(u) : \big(\eta \ot (a\vee b) \big)\hspace{1pt} &=\hspace{1pt}D^2_{ab}(\eta \cdot u)\hspace{1pt} =\hspace{1pt} D_b \big(\mathrm{G}(u): (\eta \ot a)\big),
\end{split}
\]
a.e.\ on $\Om$, where $D_a$, $D^2_{ab}$ are the usual directional derivatives. In general, the fibre spaces are strictly larger than their ``non-degenerate" counterparts and there are elements of them which are not even $W^{1,1}_{\text{loc}}$. For instance, take $\A=\eta \ot a \ot \eta \ot a $, $|a|=1$. Then, for any $f\in W^{2,2}(\R)$, $g \in L^2(\R^n)$ and $\zeta \in C^\infty_c(\R^n)$, the map $u(x) := \zeta (x)\hspace{1pt} \big[f(a\cdot x)\hspace{1pt} +\hspace{1pt}  g\big([I-a\ot a] x \big)\big]\eta$ is an element of  $\mathscr{W}^{2,2}(\Om,\Si)$ arising from this $\A$, but $D_b(\eta \cdot u)$ may not exist in $L^2$ for any $b \hspace{1pt} \bot \hspace{1pt}  a$. Similarly to the second order case, we may also define
\beq \label{4.5A}
\begin{split}
\mathscr{W}_0^{1,2}(\Om,\Si)\hspace{1pt} :=\hspace{1pt} \overline{\hspace{1pt} \text{Proj}_{ L^2\left(\Om,\Si \by \Pi \right) }\hspace{1pt}  \tilde{W}_0^{1,2}(\Om,\R^N) \hspace{1pt} }^{\|\cdot\|_{L^2(\Om)} },
\end{split}
\eeq
equipped with the obvious respective norm $\|\cdot\|_{\mathscr{W}^{1,2}(\Om)}$. Further functional properties of the fibre spaces (traces, Poincar\'e inequality) needed for the proof of Theorem \ref{theorem30} will be discussed after its statement. The fibre space $\big(\mathscr{W}^{2,2}\cap\mathscr{W}_0^{1,2}\big)(\Om,\Si)$ is the appropriate setup within which we obtain compactness and uniqueness of $\mD$-solutions for the Dirichlet problems \eqref{4.1}-\eqref{4.2}, by utilising the hypotheses introduced in the next paragraph.

\ms

\noi \textbf{Degenerate ellipticity and decomposability.} Now we introduce our ellipticity hypothesis for \eqref{4.1} and a condition for tensors $\A \in \R^{Nn \by Nn}_s$ guaranteeing their range $\Pi$ is spanned by rank-one directions.

\begin{definition}[Degenerate ellipticity]  \label{K-Condition} The Carath\'eodory map $\mF :  \Om \by \R^{Nn^2}_s \larrow \R^N$ is called degenerate elliptic when there exists a rank-one non-negative $\A \in \R^{Nn \by Nn}_s$, constants $B,C\geq 0$ with $B+C<1$ and a positive
function $A$ satisfying $A,1/A \in L^\infty(\Om)$ such that
\[
\Big|\A:\Z \hspace{1pt}- \hspace{1pt} A(x)\Big(\mF(x,\X+\Z) -\mF(x,\X)\Big) \Big|\hspace{1pt} \leq \ B \nu\hspace{1pt} |\Xi \hspace{1pt}\Z|\ +\ C\hspace{1pt}|\A :\Z|,
\]
for a.e.\ $x\in \Om$ and all $\X,\Z \in \R^{Nn^2}_s$. We moreover require $\mF$ to be valued in the subspace $\Si\sub \R^N$, i.e.\ $\mF(x,\X)\in \Si$,  for a.e.\ $x\in \Om$ and all $\X\in \R^{Nn^2}_s$.
\end{definition}

Definition \ref{K-Condition} is an extension to the degenerate realm of the strict ellipticity assumption introduced in \cite{K9}. In the elliptic case we have $\Si=\R^N$, $\Pi=\R^{Nn}$ and $\Xi=\R^{Nn^2}_s$ and then \eqref{4.1} is solvable in the class of strong solutions (for more details see \cite{K9}). If $\A_{\al i \be j}=\de_{\al \be}\de_{ij}$ and $A(x)=\text{const}$, we reduce to the classical notion introduced by Campanato (\cite{C1}-\cite{C3}). It is easy to exhibit non-trivial examples of Carath\'eodory maps satisfying Definition \ref{K-Condition}, see Remark \ref{all remarks}IV)  that follows. It is quite restrictive, but classical examples (\cite{LU}) show that even if $N=1$, extra assumption are needed for the linear uniformly elliptic equation $\A:\D^2 u =f $ if $\A$ is discontinuous. Below is the structural hypothesis for tensors.

\begin{definition}[Decomposability] \label{SH} We call $\A \in \R^{Nn \by Nn}_s$ \emph{decomposable} when it can be written as $\A_{\al i \be j} = B^1_{\al \be}A^1_{ij}\hspace{1pt} +\hspace{1pt} \cdots\hspace{1pt} +\hspace{1pt} B^N_{\al \be}A^N_{ij}$ and:

\noi i) The matrices $\{B^1,...,B^N\} \sub \R^{N^2}_s$ are non-negative with ranges $\Si^1,...,\Si^N$ mutually orthogonal in $\R^N$.

\noi ii) The matrices $\{A^1,...,A^N\} \sub \R^{n^2}_s$ are non-negative and if $\la^\ga_{i_0}$ is the smallest positive eigenvalue of $A^\ga$, the eigenspaces $N\big( A ^\ga-\la^\ga_{i_0}I \big)$ have non-trivial intersection.
\end{definition}

We discuss certain implications of these hypotheses and some examples after the main result which we give right next.

\ms

\noi \textbf{$\mD$-solutions for fully nonlinear degenerate elliptic systems.} Below is the principal result of this section followed by some relevant comments.

\begin{theorem}[Existence-Uniqueness-Partial Regularity] 
\label{theorem30} Let $\Om \sub \R^n$ be a strictly convex bounded domain with $C^2$ boundary and $\mF:\Om\by \R^{Nn^2}_s\larrow \R^N$ a map which satisfies Definition \ref{K-Condition} with respect to a decomposable $\A$ (Definition \ref{SH}). Let also $\Xi,\Pi,\Si$ be given by \eqref{4.3} and suppose that $|\mF(\cdot,0)|\in L^2(\Om)$. Then, the problem 
\[
\left\{
\begin{split}
\mF(\cdot,\D^2 u)\hspace{1pt} &=\hspace{1pt} f, \ \ \text{ in }\Om,\\
u\hspace{1pt}&=\hspace{1pt} 0, \ \ \text{ on }\p\Om, 
\end{split}
\right.
\]
has, for any $f\in L^2(\Om,\Si)$, a unique $\mD$-solution $u : \R^n \supseteq \overline{\Om} \larrow \R^N$ in the fibre space $(\mathscr{W}^{1,2}_0 \cap \mathscr{W}^{2,2})(\Om,\Si)$ (see \eqref{4.4A}, \eqref{4.5A}) with respect to certain orthonormal frames (see \eqref{2.2}) depending only on $\mF$ (Definition \ref{definition13}). In particular, $\mH^{n-1}$-a.e.\ $x\in \p\Om$ is a vanishing Lebesgue point of $u$, whilst for any $\mD ^2 u \in \mY\big(\Om,\smash{\overline{\R}}^{Nn^2}_s\big)$
\[
\ \ \ \int_{\smash{\overline{\R}}^{Nn^2}_s} \Phi(\X)\hspace{1pt} \Big(\mF(x,\X)-f(x)\Big)\hspace{1pt} d[\mD ^2 u(x)](\X)\hspace{1pt} =\hspace{1pt} 0, \quad \text{a.e.\ $x\in \Om$},
\] 
for any $\Phi\in C_c\big( {\R}^{Nn^2}_s \big)$.
\end{theorem}

\begin{remark} \label{all remarks}
I) [\textbf{Compatibility}] $f$ must be valued into $\Si$ because this is a compatibility condition due to the degeneracy of the problem. For example, the $2\by 2$ system $\De u_1=f_1$, $0 =f_2$ has no solution whatsoever in any sense unless $f_2\equiv 0$.

\smallskip

\noi II) [\textbf{Partial regularity}] The solution we obtain in Theorem \ref{theorem30} possess differentiable projections along certain rank-one lines, but in general this can not be improved further. In particular, \textit{the solution is not partially regular in the standard sense} of being more regular on a subset of the domain with full measure. For, choose any $f\in C(\overline{D})$ not weakly differentiable with respect to $x_1$ for any $x_2$ over the unit disc of $\R^2$. Then, the problem 
\[
\text{$D^2_{22}u=0$ in $D$\ \ and \  $u=0$ on $\p\D$}
\] 
has the unique explicit $\mD$-solution (which is not in $W^{1,1}_\text{loc}(\Om)$)
\[
u(x_1,x_2)\hspace{1pt} =\hspace{1pt} -v(x_1,x_2)\ +\hspace{1pt} \int_{-\infty}^{x_2}\int_{-\infty}^{t_2}f(x_1,s_2)\hspace{1pt}ds_2\hspace{1pt}dt_2,
\]
where for $(x_1,x_2)\in \overline{D}$, we set $w(x_1,x_2):=\int_{-\infty}^{x_2}\int_{-\infty}^{t_2}f(x_1,s_2)\hspace{1pt}ds_2\hspace{1pt}dt_2$ and also
\[
\begin{split}
v(x_1,x_2)\hspace{1pt} :=\ &\frac{ x_2} {2\sqrt{1-x_1^2}} \left[ w\Big(x_1, \sqrt{1-x_1^2} \Big) - w\Big(x_1, -\sqrt{1-x_1^2} \Big)\right]\\
&  +\hspace{1pt} \frac{1 } {2}\left[ w\Big(x_1, \sqrt{1-x_1^2} \Big) + w\Big(x_1, -\sqrt{1-x_1^2} \Big)\right] .
\end{split}
\]

\noi III) [\textbf{Decomposability}] Definition \ref{SH} trivialises when either $N=1$ or $n=1$ since any non-negative matrix $A \in \R^{n^2}_s$ or $B \in \R^{N^2}_s$ satisfies it. When $\max\{N,n\}\geq2$, it is non-trivial, but in view of its constructive nature it is trivial to exhibit $\A$'s satisfying it. Also, any decomposable $\A$ must be non-negative: if $Q\in \R^{Nn}$,
\[
\begin{split}
\A : Q\ot Q\hspace{1pt}= \sum_{\ga}  
\big(  (B^\ga)^{1/2}  Q  (A^\ga)^{1/2} \big)
 \hspace{1pt}\big( (B^\ga)^{1/2}  Q  (A^\ga)^{1/2} \big)\hspace{1pt} \geq \hspace{1pt} 0.
 \end{split}
\]

\noi IV) [\textbf{Examples of nonlinearities}] Fix $\A \in \R^{Nn \by Nn}_s$ and an $f\in C^{0,1}\big(\R^{Nn^2}_s, \R^N)$ with Lipschitz constant $\Lip(f)$. Then, for any positive $A$ with $A,1/A \in L^\infty(\Om)$,
\[
\mF(x,\X)\hspace{1pt} :=\hspace{1pt} \big(A(x)\big)^{-1}\Big[(1+\ga)\A :\X\hspace{1pt}+\hspace{1pt} \Si f(\Xi \X) \Big]
\]
satisfies Definition \ref{K-Condition} when $\nu |\ga|+ \Lip(f)<\nu$. Linear examples satisfying Definition \ref{K-Condition} are given by any $\A : \R^n \supseteq \Om \larrow \R^{Nn \by Nn}_s$ measurable such that $\left| \big(\A \hspace{1pt} -\hspace{1pt}  A(x)\A(x)\big):\Z \right|  \leq  B\nu\hspace{1pt}|\Xi\hspace{1pt}\Z|$ for some $0<B<1$ and $A$ positive such that $A,1/A \in L^\infty(\Om)$.

\smallskip

\noi V) [\textbf{Partial monotonicity}] Under the assumptions of Theorem \ref{theorem30}, $\mF$ satisfies
\beq \label{monotonicity}
\left\{
\begin{array}{c}
\text{For a.e.\ $x\in \Om$, $\mF(x,\cdot)$ is constant along the subspace $\Xi^\bot$, i.e.}\ms\\
  \mF(x,\X)\hspace{1pt} =\hspace{1pt} \mF(x,\Xi \hspace{1pt} \X),\ \ \ \X \in \R^{Nn^2}_s.
\end{array}
\right.
\eeq
Here $\Xi$ is as in \eqref{4.3}. Condition \eqref{monotonicity} is \emph{strictly weaker} than the decoupling condition $\mF_\al(\X) = \mF_\al(\X_\al)$ assumed in vector-valued viscosity solutions. To see it, we first note that if $\A$ satisfies Definition \ref{SH}, then $\Xi^\bot \subseteq N\big(\A : \R^{Nn^2}_s\ri \R^{N} \big)$ (this is established in the proof, see \eqref{inclusion} below). By using that $\A \!:\! \X =\A\!:\!(\Xi \hspace{1pt}\X)$, we have $\A\!:\!\Z=0$ and $\Xi\hspace{1pt}\Z=0$ when $\Z \in \Xi^\bot$. Hence, Definition \ref{K-Condition} gives $\big|A(x)\big(\mF(x,\X+\Z) -\mF(x,\X)\big) \big| \leq 0$ for all $\Z \in \Xi^\bot$ and $\X \in \R^{Nn^2}_s$. 
\end{remark}

Next we gather some properties of the fibre spaces essentially proved in \cite{K10} but without the formalism of the fibre spaces.

\begin{remark}[Basic properties of the fibre Sobolev space counterparts, cf.\ \cite{K10}] \label{fibre spaces}

\noi (I) [\textbf{Poincar\'e inequality}] For any $\Om\Subset \R^n$, unit vectors $a$, $\eta$ and $u\in W^{1,2}_0(\Om,\R^N)$, we have 
\[
\|\eta \cdot u\|_{L^2(\Om)}\hspace{1pt} \leq\hspace{1pt} \diam(\Om)\hspace{1pt} \big\|D_a(\eta \cdot u)\big\|_{L^2(\Om)}.
\]

\noi (II) [\textbf{Norm equivalence}]  The seminorm $\|G^2(\hspace{1pt} \cdot \hspace{1pt} ) \|_{L^2(\Om)}$ on the fibre space $(\mathscr{W}^{1,2}_0 \cap \mathscr{W}^{2,2})(\Om,\Si)$  (see \eqref{4.4A}, \eqref{4.5A})  is equivalent to its natural norm 
\[
\| \cdot \|_{\mathscr{W}^{2,2}(\Om)}=\hspace{1pt} \| \cdot \|_{L^2(\Om)}+ \|G(\hspace{1pt} \cdot \hspace{1pt} ) \|_{L^2(\Om)} + \|G^2(\hspace{1pt} \cdot \hspace{1pt} ) \|_{L^2(\Om)} .
\]

\noi (III) [\textbf{Trace operator}]  If $\Om\Subset \R^n$ is strictly convex and $a\in \R^n\set\{0\}$, then there is a closed set $E\sub \p\Om$ with $\mH^{n-1}(E)=0$ such that for any $\Ga\Subset \p\Om\set E$, we have
\[
\| v \|_{L^2(\Ga,\mH^{n-1})} \hspace{1pt} \leq\hspace{1pt} C \Big( \| v  \|_{L^2(\Om)}  \hspace{1pt}+\hspace{1pt} \big\| D_a v \big\|_{L^2(\Om)} \Big),
\] 
for some universal $C=C(\Ga)>0$ and all $v\in C^1(\overline{\Om})$. Hence, there is a well-defined trace operator $T: \mathscr{W}^{1,2}(\Om,\R^N)\ri L^2_{\text{loc}}(\p\Om\set E,\mH^{n-1};\R^N)$. 
\end{remark}

For the proof of Theorem \ref{theorem30} we need an important estimate established next.

\subsection{A priori degenerate estimates} \label{A priori degenerate hessian estimates} Herein we establish an a priori estimate for strong solutions in $(W^{2,2}\cap W^{1,2}_0)(\Om,\R^N)$ of a regularisation of the linear system $\A : \D^2 u=f$ when $\A$ is decomposable. This is a generalisation of the elliptic estimate of \cite{K9} (the latter extending the Miranda-Talenti inequality) to the \emph{degenerate} realm.

\begin{theorem}[Degenerate hessian estimate] \label{theorem31} Suppose $\Om\Subset \R^n$ is a convex $C^2$ domain and $\A\in \R^{Nn \by Nn}_s$ satisfies Definition \ref{SH}, $n,N\geq 1$. If $\Xi $, $\nu$ are as in \eqref{4.3}, for any $u\in (W^{2,2}\cap W^{1,2}_0)(\Om,\R^N)$ and $\e\geq 0$ we have the estimate
\[
\big\| \Xi\hspace{1pt}\D^2 u \big\|_{L^2(\Om)}\hspace{1pt} \leq\hspace{1pt} \frac{1}{\nu} \hspace{1pt} \big\| \A^{(\e)}\!: \D^2 u \big\|_{L^2(\Om)}
\]
and also the property
\beq \label{inclusion}
\Xi \hspace{1pt}  \supseteq \hspace{1pt} N\Big(\A : \R^{Nn^2}_s\ri \R^N \Big)^\bot.
\eeq
The tensor $\A^{(\e)}$ is the rank-one (strictly) positive regularisation of $\A$ given by $ \A^{(\e)}_{\al i \be j }:= \sum_{\ga=0}^N B^{(\e)\ga}_{\al \be}\hspace{1pt} A^{(\e)\ga}_{i j}$, where $B^\ga$, $A^\ga$ are as appearing in Definition \ref{SH} and
\[
\begin{split}
& B^{(\e)\ga}  \hspace{1pt} :=\hspace{1pt} \left\{
\begin{split}
& B^{\ga} , \hspace{100pt} \ga=1,...,N,\\
& \e I-\e\big(B^1+\cdots+B^N \big), \ \ \ \hspace{1pt}\ga=0,
\end{split}
\right. 
\\
& A^{(\e)\ga}  \hspace{1pt} :=\hspace{1pt} \left\{
\begin{split}
& A^{\ga} +\e I, \hspace{80pt} \ga=1,...,N,\\
& \e I, \hspace{105pt}  \ga=0.
\end{split}
\right.
\end{split}
\]
\end{theorem}

Note that in the vectorial case $N\geq 2$ of Theorem \ref{theorem31}, the ``correct" approximation is \emph{not} the vanishing viscosity one, although it reduces to such when $N=1$.

\BPT \ref{theorem31}. The first step is to prove a scalar version of the theorem.

\begin{claim} \label{claim34} Let $\Om\Subset \R^n$ be $C^2$ and convex and let also $A\geq 0$ in $\R^{n^2}_s$. Then, there exists a subspace $H  \sub \R^{n^2}_s$ such that
$H \supseteq N\left(A : \R^{n^2}_s\ri \R \right)^\bot$ and  for any $u\in (W^{2,2}\cap W^{1,2}_0)(\Om)$ and any $\e\geq 0$ we have the estimate
\[
\big\| H \D^2 u \big\|_{L^2(\Om)}\hspace{1pt} \leq\hspace{1pt} \frac{1}{\nu(A)} \hspace{1pt} \big\| A : \D^2 u \hspace{1pt}+\hspace{1pt} \e\De u\big\|_{L^2(\Om)}
\]
where $\nu(A):=\min_{|a|=1,\hspace{1pt} a\in T} \{A:a\ot a \}$ and $T := R\big(A : \R^{n}\ri \R^n \big)$.
\end{claim}

\BPC \ref{claim34}. By the Spectral theorem, we can find a diagonal matrix $\La$ with entries $0\leq \la_1\leq ... \leq \la_n$ and $O \in O(n)$ such that $A=O\La^{1/2}\hspace{1pt} (O\La^{1/2})^\top$ and
\[
\La\hspace{1pt} =\hspace{1pt} \left[
\begin{array}{c|c}
{0} & 0 \\ \hline
0 &  { 
\begin{array}{ccc}
\la_{i_0}  & & 0 \\ 
 & \ddots &  \\
0 &  & \la_{n} 
\end{array}
}
\end{array}
\right].
\]
Evidently, $\big\{ \la_1,..., \la_n\big\} = \big\{0,...,0,\la_{i_0},..., \la_n\big\}$ are the eigenvalues of $A$ and $\la_{i_0}$ is the smallest positive eigenvalue. We also fix $\e\geq 0$ and set
\beq \label{4.4}
\Th \hspace{1pt} :=\hspace{1pt} \big(\La\hspace{1pt} +\hspace{1pt} \e I\big)^{1/2}, \quad \Ga\hspace{1pt} :=\hspace{1pt} O\Th .
\eeq
Then, since $A$ equals $O \La O^\top$ and $\Th$ is symmetric, we have 
\beq  \label{4.5}
A\hspace{1pt} +\hspace{1pt} \e I\hspace{1pt} =\hspace{1pt} O \La O^\top  +\hspace{1pt} O(\e I)O^\top\hspace{1pt} =\hspace{1pt} O\Th\hspace{1pt} (O\Th)^\top \hspace{1pt} =\hspace{1pt} \Ga\Ga^\top  
\eeq
and also $\nu(A)= \la_{i_0}$ ($\nu(A)$ is defined in the statement). We define 
\beq  \label{4.6}
\begin{split}
H^0 &:=\hspace{1pt} \left\{X \in \R^{n^2}_s\ : \ X   \hspace{1pt} =\hspace{1pt} 
\left[
\begin{array}{c|c}
0 & 0  \\ \hline
0  &    (X_{ij})^{j=i_0,...,n}_{i=i_0,...,n}
\end{array}
\right]
\right\},\\
H\hspace{1pt} &:=\hspace{1pt} \Big\{X \in \R^{n^2}_s\ : \ O^\top \! X O \in H^0 \Big\}
\end{split}
\eeq
and claim the following algebraic inequality:
\beq   \label{4.7a}
 \big|\Th X \Th \big|\hspace{1pt} \geq\hspace{1pt} \nu(A)\hspace{1pt} \big|H^0 X\big|,\ \ \ X \in \R^{n^2}_s.
\eeq
Indeed, since $\Th_{ij}=0$ when $i\neq j$ and $\Th_{ii}=\sqrt{\la_i+\e}$, in view of \eqref{4.6} we have
\[
\begin{split}
 \big|\Th X \Th \big|^2 &= \sum_{i,j,k,l,p,q=1}^n \big(\Th_{ik} X_{kl} \Th_{lj} \big)\hspace{1pt} \big(\Th_{ip} X_{pq} \Th_{qj} \big) \hspace{1pt}=\\
 &=\sum_{i,j=1}^n \big(\Th_{ii} X_{ij} \Th_{jj} \big)^2\  \geq\sum_{i,j=i_0}^n \big( {\la_i +\e}\big)\hspace{1pt} ( X_{ij})^2 \big( {\la_j +\e}\big) \hspace{1pt} \geq \\
  &\geq \ (\la_{i_0})^2 \sum_{i,j=i_0}^n ( X_{ij}  )^2\ =\ \nu(A)^2\hspace{1pt} \big|H^0 X\big|^2.
\end{split}
\]
Hence, \eqref{4.7a} has been established. In order to conclude, the goal is to reduce to the classical Miranda-Talenti inequality (see \cite{M,T,K9}) which says that for $U\Subset \R^n$ convex $C^2$ domain and any $ v\in (W^{2,2}\cap W^{1,2}_0)(U)$, we have
\beq \label{4.7}
\big\| D^2v \big\|_{L^2(U)}\hspace{1pt} \leq\hspace{1pt}  \big\|  \De v\big\|_{L^2(U)}.
\eeq
It suffices to consider $\e>0$ since the case $\e=0$ follows by letting $\e\ri0$. Given any $u\in C^2(\overline{\Om})\cap C^1_0(\Om)$, we set $U:= \Ga^{-1}\Om$ and $v(x)\hspace{1pt} :=\hspace{1pt} u(\Ga x)$, $x\in U$. Then, $ \D^2_{ij}v(x)= \sum_{p,q=1}^n D^2_{p q}u(\Ga x)\hspace{1pt}\Ga_{pi} \hspace{1pt} \Ga_{qj}$ and hence, by \eqref{4.4} and \eqref{4.5} we obtain
 \beq   \label{4.8}
 \begin{split}
 D^2v(x)\ &=\ \Ga^\top \D^2 u(\Ga x)\hspace{1pt} \Ga\ =\ \Th \Big( O^\top D^2 u(\Ga x)\hspace{1pt} O \Big) \Th , \\
 \De v(x)\ &=\  \D^2 u(\Ga x) : \Ga \Ga^\top\ =\  \D^2 u(\Ga x) : (A\hspace{1pt} +\hspace{1pt} \e I).
\end{split}
\eeq
Now note that since $\Om$ is a $C^2$ bounded convex domain, $U$ is a $C^2$ bounded convex domain as well as image of such a set under a linear invertible mapping. We now apply \eqref{4.7} to $v$ over $U\sub \R^n$ and in view of \eqref{4.8} we obtain 
\[
\begin{split}
\int_U \big| \D^2 u(\Ga x) : (A\hspace{1pt} +\hspace{1pt} \e I) \big|^2d & x \hspace{1pt} \geq\hspace{1pt}  \int_U \Big|  \Th \Big( O^\top D^2 u(\Ga x)\hspace{1pt}O \Big) \Th \Big|^2dx \\
&\overset{\eqref{4.7a}}{\geq}  \nu(A)^2 \hspace{1pt} \int_U \Big|  H^0 \Big( O^\top D^2 u(\Ga x)\hspace{1pt}O \Big)   \Big|^2 dx .
\end{split}
\]
By the change of variables $y:=\Ga x$ and by using that $O$ is orthogonal, we get
\beq \label{4.9}
\big\| \D^2 u : (A\hspace{1pt} +\hspace{1pt} \e I) \big\|_{L^2(\Om)}\hspace{1pt} \geq \  \nu(A) \hspace{1pt} \Big\|  O\hspace{1pt} \Big( H^0 \left( O^\top D^2 u \hspace{1pt}O \right) \Big) \hspace{1pt} O^\top  \Big\|_{L^2(\Om)}. 
\eeq
Now we claim that the orthogonal projection on the subspace $H \sub \R^{n^2}_s$ is given by
\beq \label{4.10}
H\hspace{1pt}X\hspace{1pt} =\hspace{1pt} O\hspace{1pt} \Big( H^0 \left( O^\top X \hspace{1pt}O \right) \Big) \hspace{1pt} O^\top .
\eeq
Once \eqref{4.10} has been established, the desired estimate follows from \eqref{4.9}, \eqref{4.10}  and a standard density argument. Indeed, if $K$ denotes the linear operator defined by the right hand side of \eqref{4.10}, for any $X \in \R^{n^2}_s$ we have
\[
\begin{split}
& K\big( K\hspace{1pt} X \big) \hspace{1pt}  =\hspace{1pt} O\hspace{1pt} \Big( H^0 \left( O^\top O\hspace{1pt} \Big( H^0 \left( O^\top X \hspace{1pt} O \right) \Big) \hspace{1pt} O^\top \hspace{1pt}O \right) \Big) \hspace{1pt} O^\top =\\
& =\hspace{1pt} O\hspace{1pt} \Big( H^0  H^0 \big( O^\top X \hspace{1pt}O \big)   \Big) \hspace{1pt} O^\top =\hspace{1pt} O\hspace{1pt} \Big( H^0\big( O^\top X \hspace{1pt}O \big)  \Big) \hspace{1pt} O^\top  =\hspace{1pt} K\hspace{1pt} X.
\end{split}
\]
Hence, $K^2=K$. Moreover, $K$ is symmetric as a map $\R^{n^2}_s \larrow \R^{n^2}_s$: by using that  $H^0$ is symmetric, we have
\[
\begin{split}
&(K\hspace{1pt}X) : Y\hspace{1pt} =\hspace{1pt} \Big( O\hspace{1pt} \Big( H^0 \left( O^\top X \hspace{1pt}O \right) \Big) \hspace{1pt} O^\top \Big): Y \ =\ H^0 \left( O^\top X \hspace{1pt}O \right) :  \left( O^\top Y \hspace{1pt}O \right)=\\
 &=\hspace{1pt} \left( O^\top X \hspace{1pt}O \right) :  H^0\left( O^\top Y \hspace{1pt}O \right) \ =\ X:\Big( O\hspace{1pt} \Big( H^0 \left( O^\top Y \hspace{1pt}O \right) \Big) \hspace{1pt} O^\top \Big) \hspace{1pt} =\hspace{1pt} X:(K\hspace{1pt}Y),
\end{split}
\]
for $X,Y \in \R^{n^2}_s$. Hence, \eqref{4.10} follows. It remains to demonstrate the claimed property of $H$. To this end, fix $X \hspace{1pt} \bot \hspace{1pt} H$. Then, the projection of $X$ on $H$ vanishes and as a result of \eqref{4.10} we obtain $H^0(O^\top X\hspace{1pt} O)=0$. By recalling that $A=O\hspace{1pt} \La\hspace{1pt} O^\top$, we have $A:X=\La : (O^\top X\hspace{1pt}O)$ and since $\La \in H^0$, we conclude that $A:X=0$. Hence, we proved $H^\bot \sub N(A:\R^{n^2}_s\ri\R)$, as desired.   \qed

\ms

Next we characterise the space $H\sub \R^{n^2}_s$ of Claim \ref{claim34} in terms of  the range of $A$.

\begin{claim} \label{claim35} In the setting of Claim \ref{claim34}, we have the identity
\[
H\hspace{1pt} =\hspace{1pt} \spn[\Big\{a\vee b\ \Big| \ a,\hspace{1pt} b \in R\big(A:\R^n \ri\R^n\big) \Big\}]\hspace{1pt} =\hspace{1pt} T \vee T.
\]
\end{claim}

\BPC \ref{claim35}. We begin by observing that in view of \eqref{4.6}, we have $H=O\hspace{1pt} H^0 O^\top$ where $O\in O(n)$. Since  $H^0= \spn[ \big\{e^i \vee e^j\hspace{1pt} \big| \hspace{1pt} i,j=i_0,...,n \big\}]$, $H$ has a basis of the form $Oe^i \vee Oe^j$, $i,j=i_0,...,n$. We recall now that $A=O\La O^\top$ where $\La$ is a diagonal matrix with entries the eigenvalues $\{0,...,0,\la_{i_0},...,\la_n\}$ of $A$. We define the vectors $a^i := Oe^i\hspace{1pt}=\hspace{1pt} \big( O_{1i},...,O_{ni}\big)^\top$, $i=1,...,n$. Then, $\{a^1,...,a^n\}$ is an orthonormal frame of $\R^n$ corresponding to the columns of the matrix $A$ and is a set of eigenvectors of $A$. Since $\{a^{i_0},...,a^n\}$ correspond to the nonzero eigenvalues $\{\la_{i_0},...,\la_n\}$, the nullspace $N\big(A:\R^n \ri\R^n\big)$ is spanned by $\{a^1,...,a^{i_0-1}\}$ and hence $R\big(A:\R^n \ri\R^n\big)= \spn[\big\{ a^{i_0},...,a^n\big\}]$. Since $H$ has a basis of the form $\{a^i \vee a^j : i,j=i_0,...,n\}$, the claim follows.  \qed
 
 \smallskip

Now we begin working towards the vector case $N\geq2$. Let us first verify that $\A^{(\e)}$ is rank-one positive. Indeed, if $0<\e<1$, $\eta \in \R^N$, $a\in \R^n$, we have 
\[
\begin{split}
\A^{(\e)} : \eta \ot a \ot \eta\ot a\hspace{1pt} &\geq \min_{\de=0,...,N}\left( A^{(\e)\de}:a \ot a\right) \left[ \sum_{\ga=0}^N B^{(\e)\ga}:\eta \ot \eta \right] \hspace{1pt} \geq \\
&\geq\hspace{1pt}  \e\hspace{1pt} |a|^2  \left[ \sum_{\ga=1}^N B^{\ga} \hspace{1pt}+\hspace{1pt} \e \left(I-\sum_{\de=1}^N B^{\de} \right)\right] :\eta \ot \eta \hspace{1pt}\geq\hspace{1pt}  \e^2\hspace{1pt} |\eta|^2  |a|^2,
\end{split}
\]
as claimed. The next step is to characterise the range $\Pi$ of decomposable tensors $\A \in \R^{Nn \by Nn}_s$ in terms of the matrices $B^\ga$, $A^\ga$ composing $\A$.

\begin{claim} \label{claim36}  Let $\Pi \sub \R^{Nn}$ be the range of $\A : \R^{Nn} \larrow \R^{Nn}$ (see \eqref{4.3}). Then, $\Pi = \oplus_\ga \hspace{1pt}\big(\Si^\ga \ot T^\ga\big)$, where $\Si^\ga \sub \R^N$ and $T^\ga \sub \R^n$ are given by
\beq \label{4.11}
Si^\ga\hspace{1pt} =\hspace{1pt} R\left(B^\ga :\R^N \ri \R^N \right), \ \ \ 
T^\ga\hspace{1pt} =\hspace{1pt} R\left(A^\ga :\R^n \ri \R^n \right) 
\eeq
\end{claim}

\BPC  \ref{claim36}. We first observe that by Definition \ref{SH}, $\Si^\ga\hspace{1pt}\bot \hspace{1pt} \Si^\de$ if $\ga \neq \de$ and this implies that $\Si^\ga \ot T^\ga\hspace{1pt}\bot \hspace{1pt} \Si^\de \ot T^\de$ if $\ga \neq \de$. Let now $Q\in \R^{Nn}$. Then, $\A : Q$ is given in index form by $\sum_{\ga, \be,j}\hspace{1pt} B^\ga_{\al \be} \hspace{1pt} Q_{\be j}\hspace{1pt} A^\ga_{ji}$ which by \eqref{4.11} shows that $\Pi \sub {\oplus}_{\ga} \big(\Si^\ga \ot T^\ga\big)$. Conversely, let $R\in {\oplus}_{\ga}\big(\Si^\ga \ot T^\ga\big)$. Then, $R$ can be written as $R=  \sum_{\ga,\ka}\hspace{1pt} \big( B^\ga \eta^{\ka \ga}\big) \ot \big(A^\ga a^{\ka \ga} \big)$ for some $\eta^{\ka \ga} \in \Si^\ga$, $a^{\ka \ga} \in T^\ga$. We note that when $\ga\neq \de$, we have $\big( B^\de \ot A^\de \big) \big(\sum_{\ka}\hspace{1pt}  \eta^{\ka \ga} \ot   a^{\ka \ga}  \big)=$ because $\eta^{\ka \ga} \hspace{1pt} \bot\hspace{1pt} \Si^\de$ if $\ga\neq \de$. By defining $Q :=\sum_{\ga,\ka}\eta^{\ka \ga} \ot   a^{\ka \ga} $ it is immediate to verify that $\A:Q=R$. This establishes that $\Pi \supseteq {\oplus}_{\ga} \big(\Si^\ga \ot T^\ga\big)$, therefore completing the proof.     \qed

\ms

Next we prove an upper bound on $\nu(\A)$ in terms of $B^\ga$, $A^\ga$.

\begin{claim} \label{claim37}  Let $\nu$ be given \eqref{4.3} and $\Si^\ga$, $T^\ga$ by  \eqref{4.11}. Then, we have the estimate
\[
\nu \hspace{1pt} \leq\hspace{1pt} \left(  \min_{\ga} \hspace{1pt} \min_{\eta \in \Si^\ga,\hspace{1pt} |\eta|=1} \big\{B^\ga :\eta \ot \eta  \big\}\right)   \left(  \min_{\de} \hspace{1pt} \min_{a \in T^\de,\hspace{1pt} |a|=1} \big\{A^\de :a \ot a  \big\}\right).
\]
\end{claim}

\BPC  \ref{claim37}. We begin by noting that on top of the decomposability we may further assume that all the matrices $A^\ga$ have the same smallest positive eigenvalue $\la^\ga_{i_0}$ equal to $1$ for all $\ga=1,...,N$ which is realised at a common eigenvector $\bar{a}\in \R^n$. Indeed, existence of $\bar{a}$ follows from Definition \ref{SH} since the eigenspaces 
$N\big(A^\ga-\la^\ga_{i_0}I \big)$ intersect for all $\ga$ at least along a common line in $\R^n$. Further, by replacing $\{B^1,...,B^N \}$,  $\{A^1,...,A^N \}$ by the rescaled families $\{\tilde{B}^1,...,\tilde{B}^N \}$,  $\{\tilde{A}^1,...,\tilde{A}^N \}$ where $\tilde{B}^\ga := \la^\ga_{i_0} B^\ga$, $\tilde{A}^\ga := (1/\la_{i_0}^\ga) A^\ga$, we have that the new families have the same properties as the original and in addition all the new $A^\ga$ matrices have the same minimum positive eigenvalue normalised to $1$. Hence, we may assume 
\beq \label{4.12}
\ \exists\ \hspace{1pt} \bar{a}\hspace{1pt} \in \hspace{1pt} \p \mB^n_1\overset{N}{\underset{\ga=1}{\cap}}\ T^\ga \ \ :\ \ \quad \la^\ga_{i_0} =  \min_{a \in T^\ga,\hspace{1pt} |a|=1} \big\{A^\ga :a \ot a  \big\} = A^\ga : \bar{a} \ot \bar{a}=1,
\eeq
for $\ga=1,...,N$. By using \eqref{4.12}, Claim \ref{claim36} and that $\cup_\ga \big(\Si^\ga \ot T^\ga\big) \sub \oplus_\ga \big(\Si^\ga \ot T^\ga\big) $
\[
\begin{split}
\nu \hspace{1pt} &=\hspace{1pt} \min_{|\eta|=|a|=1,\hspace{1pt} \eta \ot a \in \Pi} \sum_{\de} \Big(B^\de :\eta \ot \eta \Big) \Big(A^\de :a \ot a \Big)\\
 &\leq\hspace{1pt} \min_{|\eta|=|a|=1,\hspace{1pt} \eta \ot a \in \cup_\ga (\Si^\ga \ot T^\ga)} \sum_{\de} \Big(B^\de :\eta \ot \eta \Big) \Big(A^\de :a \ot a \Big),
\end{split}
\] 
and hence
\[
\begin{split}
\nu \hspace{1pt}  &\leq\hspace{1pt} \min_\ga  \left( \min_{|\eta|=|a|=1,\hspace{1pt} \eta \ot a \in \Si^\ga \ot T^\ga} \sum_{\de} \Big(B^\de :\eta \ot \eta \Big) \Big(A^\de :a \ot a \Big) \right)
\\
&\leq\hspace{1pt} \min_\ga  \left( \min_{|\eta|=|a|=1,\hspace{1pt} \eta \ot a \in \Si^\ga \ot T^\ga} \sum_{\de} \Big(B^\de :\eta \ot \eta \Big) \Big(A^\de :a \ot a \Big) \right)\\
&\leq\hspace{1pt} \min_\ga  \left( \min_{|\eta|=1,\hspace{1pt} \eta \in \Si^\ga} \sum_{\de} \Big(B^\de :\eta \ot \eta \Big) \Big( A^\de : \bar{a} \ot \bar{a} \Big) \right)\\
&=\hspace{1pt} \min_\ga   \min_{|\eta|=1,\hspace{1pt} \eta \in \Si^\ga} \sum_{\de} \Big(B^\de :\eta \ot \eta \Big) .
\end{split}
\]
Since $B^\de :\eta \ot \eta =0$ if $\eta \in \Si^\ga$ for $\ga \neq \de$, by using \eqref{4.12} again we conclude that $\nu \leq  \min_{\ga} \min_{\eta \in \Si^\ga,\hspace{1pt} |\eta|=1} \big\{B^\ga :\eta \ot \eta  \big\}$, as desired.        \qed

\ms

Now we complete the proof by using the previous claims. We define 
\beq \label{4.13A}
\Xi\hspace{1pt} :=\hspace{1pt} \underset{\ga}\oplus\hspace{1pt} \Big( \Si^\ga \ot T^\ga \vee T^\ga \Big)\ \ \ \sub \ \R^{Nn^2}_s,
\eeq
and for brevity we set $\Xi^\ga := T^\ga \vee T^\ga$ where $\Si^\ga$, $T^\ga$ are as in \eqref{4.11}. Fix $u\in C^2(\overline{\Om},\R^N)\cap C^1_0(\Om,\R^N)$. Then, for $\ga,\al=1,...,N$, by Claims \ref{claim34}, \ref{claim35} applied to the scalar function $(\Si^\ga u)_{\al} \in C^2(\overline{\Om} )\cap C^1_0(\Om)$, we have
\[
\int_\Om \Big| \Xi^\ga D^2(\Si^\ga u)_{\al} \Big|^2\hspace{1pt} \leq\hspace{1pt} \int_\Om \Big| A^{(\e)\ga}:D^2 (\Si^\ga u)_{\al} \Big|^2,
\]
where we have used that $A^{(\e)\ga}=A^\ga+\e I$ (by the definition of $\A^{(\e)}$) and we have employed the normalisation of \eqref{4.12} which forces $\la_{i_0}^\ga = \nu(A^\ga)=1$. By summing in $\al,\ga$, the above estimate and \eqref{4.13A} give
\beq \label{4.13}
\begin{split}
\int_\Om \big| \Xi\hspace{1pt} D^2 u \big|^2\hspace{1pt} =\hspace{1pt} \int_\Om \sum_\ga \Big| \Si^\ga \ot\Xi^\ga : \D^2 u\Big|^2\ \leq\hspace{1pt} \int_\Om  \sum_\ga \Big| \Si^\ga \big( D^2 u : A^{(\e)\ga}\big)\Big|^2.
\end{split}
\eeq
We also set $C^{(\e)\ga}:=\Si^\ga \big( D^2 u : A^{(\e)\ga}\big)$ for $\ga=1,...,N$. Then, \eqref{4.13} says
\beq \label{4.14}
\int_\Om \big| \Xi\hspace{1pt} D^2 u \big|^2\hspace{1pt} \leq\hspace{1pt}  \int_\Om \sum_{\ga=1}^N  \big| C^{(\e)\ga}  \big|^2.
\eeq
By the definition of $\A^{(\e)}$, we have that $B^{(\e)\ga} \hspace{1pt} \bot \hspace{1pt} B^{(\e)\de}$ for $\ga \neq \de$ in $\{0,1,...,N\}$. By using this fact, we calculate
\[
\begin{split}
\big| \A^{(\e)} : \D^2 u \big|^2 
&=\hspace{1pt} \left( \sum_{\ga=0}^N  B^{(\e)\ga} \hspace{1pt} \big( D^2 u : A^{(\e)\ga} \big) \right) \cdot \left( \sum_{\de=0}^N  B^{(\e)\de} \hspace{1pt} \big( D^2 u : A^{(\e)\de} \big) \right)\\
 &=\hspace{1pt} \sum_{\ga=0}^N  \Big(  B^{(\e)\ga} \hspace{1pt} \big( D^2 u : A^{(\e)\ga}\big) \Big)\cdot\Big( B^{(\e)\ga}  \big( D^2 u : A^{(\e)\ga}\big) \Big)
 \end{split}
\]
and hence
\[
\begin{split}
\big| \A^{(\e)} : \D^2 u \big|^2 \hspace{1pt} &=\hspace{1pt}  \Big|  B^{(\e)0} \hspace{1pt} \big( D^2 u : A^{(\e)0}\big) \Big|^2\hspace{1pt}+\hspace{1pt} \sum_{\ga=1}^N  \Big|  B^{(\e)\ga} \hspace{1pt} \big( D^2 u : A^{(\e)\ga}\big)
 \Big|^2  \geq
 \end{split}
\]

 \[
\begin{split}
 &\geq \sum_{\ga=1}^N  \Big|  B^{\ga} \hspace{1pt} \big( D^2 u : A^{(\e)\ga}\big)
 \Big|^2 \hspace{1pt} =\ \sum_{\ga=1}^N  \big|  B^{\ga} \hspace{1pt} C^{(\e)\ga} \big|^2 \hspace{1pt} \geq \hspace{1pt} \sum_{\ga=1}^N  \hspace{1pt} \max_{|\eta|=1} \Big(  B^{\ga} : \big(C^{(\e)\ga} \ot \eta \big)\Big)^2\geq	
 \\
 &\geq\hspace{1pt} \sum_{\ga=1}^N  \hspace{1pt} \Big(  B^{\ga} : \left(\sgn(C^{(\e)\ga}) \ot \sgn(C^{(\e)\ga})  \right) \Big)^2 \big| C^{(\e)\ga}  \big|^2  .
\end{split}
\]
As a result, we obtain
\beq \label{4.16}
\begin{split}
\big| \A^{(\e)} : \D^2 u \big|^2 \hspace{1pt} &\geq\hspace{1pt}   \left( \min_{\de=1,...,N} \min_{|\eta|=1,\hspace{1pt}\eta \in \Si^\de} \big\{ B^{\de} : \eta\ot \eta\big\} \right)^2 \sum_{\ga=1}^N  \big| C^{(\e)\ga}  \big|^2.
\end{split}
\eeq
By using Claim \ref{claim37} (and also the normalisation condition \eqref{4.12}), \eqref{4.16} gives
\beq \label{4.17}
\int_\Om \big| \A^{(\e)} : \D^2 u \big|^2 \geq\hspace{1pt} \nu^2 \int_\Om \hspace{1pt} \sum_{\de=1}^N \hspace{1pt} \big| C^{(\e)\de}  \big|^2.
\eeq
Hence, by \eqref{4.17} and \eqref{4.14} we obtain the desired estimate for smooth $u$, the general case following by a density argument. We complete the proof by showing that $\Xi$ satisfies \eqref{inclusion}. If $\X \hspace{1pt} \bot \hspace{1pt} \Xi $, by \eqref{4.13A} $\X$ is normal to $\Si^\ga \ot H^\ga$ for any $\ga=1,...,N$, where $H^\ga:= T^\ga \vee T^\ga$. Hence the projection of $\X$ on $\Si^\ga \ot H^\ga$ vanishes: $(\Si^\ga \ot H^\ga)\X=0$. By Claim \ref{claim34} we have $A^\ga :X = A^\ga : (H^\ga X)$ if $X\in \R^{n^2}$. Hence, $B^\ga (\X : A^\ga)=0$ for $\ga=1,...,N$ and by summing in $\ga$ we obtain $\A :\X=0$. Thus, we have shown $\Xi^\bot  \sub N\big(\A : \R^{Nn^2}_s\ri \R^N \big)$ as desired.            \qed

\ms

\subsection{Proof of the main result} \label{Proof of the main result} Now we establish our second main result by utilising the a priori estimate of Subsection \ref{A priori degenerate hessian estimates}. 

\BPT \ref{theorem30}. The fist step is to prove existence of a map in the fibre space $(\mathscr{W}^{2,2}\cap\mathscr{W}^{1,2}_0)(\Om,\Si)$ solving in a certain sense the linear problem.

\begin{claim} \label{claim41} In the setting of Theorem \ref{theorem30} and under the same assumptions, for any $f\in L^2(\Om,\Si)$ there exists a unique $u\in (\mathscr{W}^{2,2}\cap\mathscr{W}^{1,2}_0)(\Om,\Si)$ such that $\A : \mathrm{G}^2(u) = f$ a.e.\ on $\Om$, where $\mathrm{G}^2(u)$ is the fibre hessian of $u$.
\end{claim}

\BPC \ref{claim41}. The proof is based on the approximation by strictly elliptic systems and relies on  the stable estimate of Theorem \ref{theorem31}. Let $\A^{(\e)}$ be the approximation of $\A$ of Theorem \ref{theorem31} and consider for a fixed $f\in L^2(\Om,\Si)$ the system $\A^{(\e)} : \D^2 u^\e = f$, a.e. on $\Om$. By standard lower semicontinuity and regularity results (see e.g.\ \cite{D, GM}), the problem has for any $\e>0$ a unique strong a.e.\ solution $u^\e \in (W^{2,2}\cap W^{1,2}_0)(\Om,\R^N)$. By  Theorem \ref{theorem31} and Remark \ref{fibre spaces}, we have the estimate
\[
\big\|\Si \hspace{1pt} u^\e \big\|_{L^2(\Om)}\hspace{1pt}+\hspace{1pt} \big\|\Pi\hspace{1pt} D u^\e \big\|_{L^2(\Om)}\hspace{1pt}+\hspace{1pt} \big\|\Xi\hspace{1pt} \D^2 u^\e \big\|_{L^2(\Om)} \hspace{1pt} \leq\hspace{1pt} \frac{C}{\nu} \hspace{1pt} \big\| f \big\|_{L^2(\Om)}
\]
for some universal $C>0$. By the definition of $ (\mathscr{W}^{2,2}\cap \mathscr{W}^{1,2}_0)(\Om,\Si)$ (\eqref{4.4A},\eqref{4.5A}), there exists $u$ such that $\big(\Si \hspace{1pt} u^\e ,\Pi\hspace{1pt} D u^\e ,\Xi\hspace{1pt} D^2 u^\e \big)\weak \big(u,\mathrm{G}(u),\mathrm{G}^2(u)\big)$, along a sequence $\e_k\ri 0$ in $L^2$. Now we pass to the weak limit in the equations. By the form of the approximation $\A^{(\e)}$ and Definition \ref{SH}, we have
\[
\sum_{\ga=1}^N \hspace{1pt} B^{(\e)\ga} \left( \D^2 u^\e : A ^{(\e)\ga} \right)\hspace{1pt} =\hspace{1pt} f\hspace{1pt}-\hspace{1pt} B^{(\e)0} \big( \D^2 u^\e : A ^{(\e)0}\big),
\]
a.e.\ on $\Om$. By using that $B^{(\e)\ga}=B^\ga$ for $\ga=1,...,N$ and that $B^{(\e)0} \hspace{1pt} \bot \hspace{1pt} B^1+\cdots+B^N$, we may project the system above on  the range of $ B^1+\cdots+B^N$ which we denote by $\Si$. Then, since $\Si f=f$ and $A ^{(\e)\ga}=A ^{\ga}+\e I$, we obtain
\[
\sum_{\ga=1}^N \hspace{1pt} B^{\ga} \Big(  \e\De u^\e \ +\ \D^2 u^\e : A ^{\ga}\Big)\hspace{1pt} =\hspace{1pt} f,
\]
a.e.\ on $\Om$. Moreover, by \eqref{inclusion} (and in view of Remark \ref{all remarks}), we deduce that for any $\phi \in C^\infty_c(\Om,\R^N)$, integration by parts gives
\[
\int_\Om   \Big( \A : \big( \Xi\hspace{1pt} \D^2 u^\e \big) \hspace{1pt} -\hspace{1pt}f\Big)\cdot \phi \hspace{1pt} =\hspace{1pt} -\hspace{1pt}\e\int_\Om \sum_{\ga=1}^N \hspace{1pt} B^{\ga} (\Si u^\e) \cdot \De \phi.
\]
By letting $\e_k\ri 0$, we obtain $\A :\mathrm{G}^2(u)=f$, a.e.\ on $\Om$. We finally show uniqueness. Let $v,w \in (\mathscr{W}^{2,2}\cap\mathscr{W}^{1,2}_0)(\Om,\Si)$ be two solutions of the system. Then, there are sequences $(v^m)_1^\infty,(w^m)_1^\infty \sub (W^{2,2}\cap W^{1,2}_0)(\Om,\R^N)$ such that  $ v^m-w^m \larrow v-w$ with respect to $\|\cdot\|_{\mathscr{W}^{2,2}(\Om)}$ as $m\ri \infty$. By assumption we have $\A : G^2(v-w)=0$ a.e.\ on $\Om$, and hence $\A : D^2(v^m-w^m)= :f^m$ a.e.\ on $\Om$ and $f^m  \ri 0$ in $L^2(\Om,\R^N)$ as $m\ri \infty$. Hence, by Theorem \ref{theorem31} and Remark \ref{fibre spaces},
\[
\begin{split}
 \|f^m  \|_{L^2(\Om)} \geq\hspace{1pt} \nu \big\|\Xi: D^2(v^m -w^m) \big\|_{L^2(\Om)} \geq\hspace{1pt}  C\big\|\Si (v^m -w^m) \big\|_{L^2(\Om)}
\end{split}
\]
and by letting $m\ri \infty$ we see that $v\equiv w$, hence uniqueness ensues.  \qed

\ms

An essential ingredient in order to pass to the nonlinear problem is the next result of Campanato (\cite{C3}, \cite{K9}) which we recall for the convenience of the reader.

\begin{lemma}[Campanato's bijectivity of near operators] \label{Campanato} Let $\mathfrak{X}\neq \emptyset$ be a set, $(X,\|\cdot\|)$ a Banach space and $\mathscr{F},\mathscr{A} : \mathfrak{X} \larrow X$ maps such that
\[
\Big\| \mathscr{F}(u)-\mathscr{F}(v)\hspace{1pt} -\hspace{1pt} \big( \mathscr{A}(u)-\mathscr{A}(v)\big)\Big\|\hspace{1pt}\leq\hspace{1pt} K\hspace{1pt} \big\|\mathscr{A}(u)-\mathscr{A}(v) \big\| 
\]
for some $K\in (0,1)$ and all $u,v \in \mathfrak{X} $. Then, if $\mathscr{A}$ is bijective, $\mathscr{F}$ is bijective as well.
\end{lemma}

Now we employ Lemma \ref{Campanato} in order to show existence of  a map in the fibre space $(\mathscr{W}^{2,2}\cap\mathscr{W}^{1,2}_0)(\Om,\Si)$ solving in a certain sense the nonlinear problem.

\begin{claim} \label{claim43} In the setting of Theorem \ref{theorem30} and under the same assumptions, for any $f\in L^2(\Om,\Si)$ there exists a unique $u\in (\mathscr{W}^{2,2}\cap\mathscr{W}^{1,2}_0)(\Om,\Si)$ such that $\mF\big(\cdot, \mathrm{G}^2(u)\big)= f$ a.e.\ on $\Om$ where $\mathrm{G}^2(u)$ is the fibre hessian of $u$.
\end{claim}

\BPC \ref{claim43}. For any fixed $u\in (\mathscr{W}^{2,2}\cap\mathscr{W}^{1,2}_0)(\Om,\Si)$, we have that $\A:\mathrm{G}^2(u)$ is in $L^2(\Om,\Si)$ because $\mathrm{G}^2(u) \in L^2(\Om,\Xi)$ and also $\A :\X$ lies is in $\Si\sub \R^N$ for any $\X \in \Xi \sub \R^{Nn^2}_s$. Moreover, by Definition \ref{K-Condition}  we have
\[
\big| \mF\big(\cdot, \mathrm{G}^2(u)\big) \big|\hspace{1pt} \leq\hspace{1pt} \Bigg( \frac{ (C+1)|\A|\hspace{1pt} +\hspace{1pt} B\hspace{1pt} \nu
 }{ {\ess\hspace{1pt}\inf}_{x\in \Om} [A(x)] }\Bigg) \hspace{1pt} |\mathrm{G}^2(u)|\hspace{1pt} +\hspace{1pt} \big| \mF\big(\cdot, 0\big) \big|,
\]
a.e.\ on $\Om$. Hence, $\mF\big(\cdot, \mathrm{G}^2(u)\big)$ is in $L^2(\Om,\Si)$ as well. The previous considerations imply that the maps
\[
\begin{split}
&\mathscr{A}\ : \ (\mathscr{W}^{2,2}\cap\mathscr{W}^{1,2}_0)(\Om,\Si) \larrow L^2(\Om,\Si), \ \ \ \mathscr{A}(u)\hspace{1pt} :=\hspace{1pt} \A :\mathrm{G}^2(u),\\
&\mathscr{F}\ : \ (\mathscr{W}^{2,2}\cap\mathscr{W}^{1,2}_0)(\Om,\Si) \larrow L^2(\Om,\Si), \ \ \  \mathscr{F}(u)\hspace{1pt} :=\mF\big(\cdot, \mathrm{G}^2(u)\big),
 \end{split}
\]
are well defined. By Claim \ref{claim41}, $\mathscr{A}$ is bijective. We complete the claim by showing that $\mathscr{F}$ is near $\mathscr{A}$ in the sense of Lemma \ref{Campanato}. For any $u,v\in (\mathscr{W}^{2,2}\cap\mathscr{W}^{1,2}_0)(\Om,\Si)$, by Definition \ref{K-Condition} and Theorem \ref{theorem31} we have
\[
\begin{split}
\Big\|A(\cdot)\hspace{1pt}& \Big( \mF\big(\cdot, \mathrm{G}^2(u)\big) - \mF\big(\cdot, \mathrm{G}^2(v)\big)\Big)\hspace{1pt} -\hspace{1pt} \A:\big(\mathrm{G}^2(u)-\mathrm{G}^2(v) \big) \Big\|_{L^2(\Om)}  \\
&\leq\hspace{1pt} B\hspace{1pt}\nu \hspace{1pt}\big\| \mathrm{G}^2(u)-\mathrm{G}^2(v) \big\|_{L^2(\Om)}  \hspace{1pt} +\hspace{1pt} C\hspace{1pt}  \big\| \A:\big(\mathrm{G}^2(u)-\mathrm{G}^2(v) \big) \big\|_{L^2(\Om)}  \\
&\leq\hspace{1pt} (B+C) \hspace{1pt}\big\| \A:\big(\mathrm{G}^2(u)-\mathrm{G}^2(v) \big) \big\|_{L^2(\Om)}.  
\end{split}
\]
Hence, $\hat{\mathscr{F}}(u):=A(\cdot)\hspace{1pt}\mF\big(\cdot, \mathrm{G}^2(u)\big)$ is bijective and since $A,1/A \in L^\infty(\Om)$, the same is true for $\mathscr{F}$. The claim ensues.          \qed

\ms

The next claim completes the proof of Theorem \ref{theorem30}.

\begin{claim} \label{claim44} In the setting of Claim \ref{claim43} and under the same assumptions, there exists an orthonormal frame $\{E^1,...,E^N\} \sub \R^N$ and for each $\al=1,...,N$ there is an orthonormal frame $\{E^{(\al)1},...,E^{(\al)n}\} \sub \R^n$ (both depending only on $\mF$) such that,  the map $u$ corresponding to $f\in L^2(\Om,\Si)$ is the unique $\mD$-solution of
$\mF(\cdot, \D^2 ) =f$ in the fibre space $(\mathscr{W}^{2,2}\cap\mathscr{W}^{1,2}_0)(\Om,\Si)$.
\end{claim}

\BPC \ref{claim44}.  \textbf{Step 1 (The frames).} By \eqref{4.3} and \eqref{4.11} there is a frame $\{E^\al | \al \}$ of $\R^N$ and for each $\al$ there is a frame $\{E^{(\al)i}| i\}$ of $\R^n$ such that each of the mutually orthogonal subspaces $\Si^\ga \sub \R^N$ is spanned by a subset of vectors $E^\al$ and for the same index $\ga$, $T^\ga$ is spanned by $\{E^{(\al)i_0},...,E^{(\al)n}\}$ which is a set of eigenvectors of $A^\ga$. By  \eqref{4.3} and \eqref{4.13A}  there are also induced frames of $\R^{Nn}$ and $\R^{Nn^2}_s$ of matrices as in \eqref{2.2}. These frames are such that a subset of the $E^{\al i j}$'s spans the subspace $\Xi$ and the rest are orthogonal to $\Xi$. 
\smallskip

\noi  \textbf{Step 2 (Sufficiency).}  Let now $u\in (\mathscr{W}^{2,2}\cap\mathscr{W}^{1,2}_0)(\Om,\Si)$ be the map of Claim \ref{claim43} which satisfies $\mF(\cdot,\mathrm{G}^2(u))=f$ a.e.\ on $\Om$. Let also us fix any infinitesimal sequence $(h_{\underline{m}})_{\underline{m}\in \N^2}$ with respect to the frames of Step 1 (see Definition \ref{definition6}) and let $\mD ^2 u$ be any diffuse hessian of $u$ arising from this sequence $\de_{\D^{2,h_{\underline{m} }}u} \weakstar \hspace{1pt} \mD^2 u$ in $\mY\big( \Om,\smash{\overline{\R}}^{Nn^2}_s \big)$  as $\underline{m}\ri \infty$, perhaps along subsequences. By the characterisation of the fibre hessian $\mathrm{G}^2(u)\in L^2(\Om,\Xi)$ in terms of directional derivatives of projections (Subsection \ref{4.2}),
\beq \label{eq}
\ \ \ \mathrm{G}^2(u)\hspace{1pt} =\sum_{\al,i,j\hspace{1pt} :\hspace{1pt} E^{\al i j} \in \Xi}  \Big(\mathrm{G}^2(u):E^{\al i j}\Big) E^{\al i j}, \ \ \ \text{ a.e.\ on }\Om,
\eeq
because the projection of $\mathrm{G}^2(u)$ along $E^{\al i j}$ is non-zero only for those  $E^{\al i j}$ spanning $\Xi$. Since $\mF$ is a Carath\'eodory map and $\mF\big(x,\mathrm{G}^2(u)(x)\big)=f(x)$ for a.e.\ $x\in\Om$, by \eqref{eq} and in view of \eqref{2.3b} we get
\[
\mF\Bigg( x, \sum_{\al,i,j\hspace{1pt} :\hspace{1pt} E^{\al i j} \in \Xi} \left[ \D^{2, h_{ m_1^2}h_{m_2^2 }}_{E^{(\al)i} E^{(\al)j}}\big( E^\al \cdot u\big) \right] (x) E^{\al i j}\Bigg)  \larrow f(x),
\]
for a.e.\ $x\in \Om$ as $\underline{m} \ri \infty$. By Remark \ref{all remarks}V), the above is equivalent to
\[
\mF\Big(x, \D^{2,h_{\underline{m}}}u(x)\Big)\hspace{1pt}=\hspace{1pt} \mF\Bigg(x, \sum_{\al,i,j } \left[ \D^{2, h_{ m_1^2}h_{m_2^2 }}_{E^{(\al)i} E^{(\al)j}}\big( E^\al \cdot u\big) \right](x) E^{\al i j}\Bigg) \larrow f(x),
\]
for a.e.\ $x\in \Om$, as $\underline{m} \ri \infty$. We set $f^{ \underline{m} }(x) := \mF\big(x, \D^{2,h_{\underline{m}}}u(x)\big) -\hspace{1pt} f(x)$ and note that we have $f^{ \underline{m} } \larrow 0$, a.e.\ on $\Om$ as $\underline{m} \ri \infty$. By the above, for any $\Phi \in C_c\big({\R}^{Nn^2}_s\big)$,
\[
\int_{\smash{\overline{\R}}^{Nn^2}_s} \Phi(\X)\hspace{1pt} \Big[ \mF(x,\X) - \big(f(x)+ f^{ \underline{m} }(x)\big) \Big] 
\hspace{1pt}d\big[ \de_{\D^{2,s_{ \underline{m} }}u(x)}\big](\X)\hspace{1pt} =\hspace{1pt} 0, \quad \text{ a.e. }x\in \Om.
\]
Since $f^{ \underline{m} } \ri0$ a.e.\ on $\Om$ as $\underline{m} \ri \infty$, we apply the Convergence Lemma \ref{lemma16A} to obtain
\[
\int_{\smash{\overline{\R}}^{Nn^2}_s} \Phi(\X)\hspace{1pt} \big[ \mF(x,\X) - f(x) \big] 
 \hspace{1pt}d\big[ \mD^2 u(x)\big](\X)\hspace{1pt} =\hspace{1pt} 0, \quad \text{ a.e. }x\in \Om,
\]
for any $\Phi \in C_c\big({\R}^{Nn^2}_s\big)$. Hence, the map $u$ of Claim \ref{claim43} is a $\mD$-solution of \eqref{4.1}.

\smallskip

\noi  \textbf{Step 3 (Necessity).} We now finish the proof by showing that any $\mD$-solution $w$ of \eqref{4.1} with respect to the frames of Step 1 which lies in the fibre space $(\mathscr{W}^{2,2}\cap\mathscr{W}^{1,2}_0)(\Om,\Si)$ actually coincides with the map $u$ of Claim \ref{claim43}. By Theorem \ref{proposition15}, we have that the $\mD$-solution $w$ can be characterised by the property that for any $R>0$, the cut off associated to $\mF$ (see Definition \ref{definition24}) satisfies $\mF\big(\cdot, \big[\D^{2,h_{\underline{m}}}w\big]^R\big) \larrow f$, a.e.\ on $\Om$ as $\underline{m}\ri \infty$. By using Remark \ref{all remarks}V), we have for any $R>0$ that
$\mF\big(\cdot, \big[\Xi\hspace{1pt} \D^{2,h_{\underline{m}}}w\big]^R\big) \larrow f$, a.e.\ on $\Om$ as $\underline{m}\ri \infty$. Since $w$ is in $(\mathscr{W}^{2,2}\cap\mathscr{W}^{1,2}_0)(\Om,\Si)$, by the properties of the fibre space we get $\Xi \hspace{1pt} \D^{2,h_{\underline{m}}}w \larrow G^2(w)$ in $L^2$ and hence a.e.\ on $\Om$ along perhaps subsequences. By passing to the limit as $\underline{m}\ri \infty$ and then as $R\ri \infty$, we obtain that $\mF(\cdot,G^2(w))=f$, a.e.\ on $\Om$. Hence, $w\equiv u$.    \qed

\smallskip

By recalling Remark \ref{fibre spaces} regarding the boundary trace values of maps in the fibre space, we conclude that the proof of Theorem \ref{theorem30} is now complete.             \qed

\begin{remark}[Regularity of $\mD$-solutions] In a sense, Claim \ref{claim44} says that all diffuse hessians of the $\mD$-solution $u$ when restricted to the subspace of non-degeneracies have the ``functional" representation $\mathrm{G}^2(u)$ \emph{inside the coefficients}. Indeed, by decomposing $\R^{Nn^2}_s=\Xi \oplus \Xi^\bot$, the restriction of any $\mD ^2 u \in \mY\big(\Om, \smash{\overline{\R}}^{Nn^2}_s\big)$ to $\Xi$ is given by the fibre hessian: $ \mD ^2 u(x)\hspace{1pt}\LL \hspace{1pt} \Xi =\de_{\hspace{1pt}G^2u(x)}$, a.e. $x\in\Om$.  Although such a simple representation is not possible in general (compare e.g.\ with Theorems \ref{theorem20}, \ref{lemma31}), it is expected that weaker versions of such results should be true (see Proposition \ref{lemma13}). 
\end{remark}

\noi \textbf{Acknowledgement.} I wish to thank Craig Evans, Robert Jensen, Charles Fefferman, Luigi Ambrosio, Roger Moser and Jan Kristensen most warmly for our inspiring scientific discussions and their appreciation of this work. I am also indebted to Alessio Figalli, Martin Hairer, Nicholas Barron, Tristan Pryer and Filippo Cagnetti for their various comments which improved the material of this paper.

\bibliographystyle{amsplain}

\end{document}